\documentclass[%
 aip,
 pof,%
 amsmath,amssymb,
preprint,%
]{revtex4-2}

\usepackage{graphicx}
\usepackage{dcolumn}
\usepackage{bm}

\usepackage[utf8]{inputenc}
\usepackage{graphicx}
\usepackage{amsmath, amsfonts, amssymb}
\usepackage{cancel}
\usepackage{physics}
\usepackage{hyperref}
\usepackage[english]{babel}
\usepackage{multirow}

\renewcommand{\vec}[1]{\boldsymbol{#1}}
\renewcommand{\L}{\mathcal{L}}

\usepackage[dvipsnames]{xcolor}

\graphicspath{{../figures/}}

\begin{document}

\preprint{}

\title{Refined radial basis function-generated finite difference analysis of non-Newtonian natural convection}%

\author{Miha Rot}
\email{miha.rot@ijs.si}
 \affiliation{%
 	Jozef Stefan International Postgraduate School, Jamova 
 	cesta 39, 1000 Ljubljana, Slovenia
 }%
\affiliation{%
	Jožef Stefan Institute, Parallel and Distributed Systems 
	Laboratory, Jamova cesta 39, 1000 Ljubljana, Slovenia
}%
\author{Gregor Kosec}%
\affiliation{%
Jožef Stefan Institute, Parallel and Distributed Systems 
Laboratory, Jamova cesta 39, 1000 Ljubljana, Slovenia
}%

\date{\today}

\begin{abstract}
In this paper we present a refined Radial Basis Function-generated Finite Difference (RBF-FD) solution for a non-Newtonian fluid in a closed differentially heated cavity. The non-Newtonian behaviour is modelled with the Ostwald-de Waele power law and the buoyancy with the Boussinesq approximation. The problem domain is discretised with scattered nodes without any requirement for a topological relation between them. This allows a trivial generalisation of the solution procedure to complex irregular three dimensional (3D) domains, which is also demonstrated by solving the problem in a two dimensional (2D) and 3D geometry mimicking a porous filter. The results in 2D are compared with two reference solutions that use the Finite volume method in a conjunction with two different stabilisation techniques, where we achieved good agreement with the reference data. The refinement is implemented on top of a dedicated meshless node positioning algorithm using piecewise linear node density function that ensures sufficient node density in the centre of the domain while maximising the node density in a boundary layer where the most intense dynamic is expected. The results show that with a refined approach, more than 5 times fewer nodes are required to obtain the results with the same accuracy compared to the regular discretisation. The paper also discusses the convergence with refined discretisation for different scenarios for up to $2 \cdot 10^5$ nodes, the impact of method parametres, the behaviour of the flow in the boundary layer, the behaviour of the viscosity and the geometric flexibility of the proposed solution procedure.
\end{abstract}

\maketitle

\section{Introduction}
Natural convection, a type of flow driven by the temperature-dependent density of a fluid, is a cornerstone of many natural and industrial processes. The most vivid examples in nature where it plays a crucial role are probably weather systems, e.g. sea and land breezes. In industry, it is of paramount importance in metal casting, various heating systems, food processing, etc.~\protect\cite{rahimi2018comprehensive} The minimal model describing natural convection involves coupled heat transfer and fluid dynamics~\protect\cite{bejan2013convection}, where simple fluids such as air and water are usually modelled as Newtonian fluids, i.e. the viscosity of the fluid is constant. However, such an approximation becomes insufficient when dealing with more complex fluids~\cite{Chhabra2010} like melts of large and complicated molecules, polymers~\cite{Wang2010}, suspensions, foams, biological fluids such as blood~\cite{Charm1965}, food~\cite{Gratao2007, WeltiChanes2005}, etc. In such fluids, also called non-Newtonian fluids, the relationship between stress and strain is no longer linear and the viscous stress becomes related to the shear rate~\cite{bingham1917investigation}. The behaviour of purely viscous non-Newtonian fluids can generally be divided into two groups based on how their viscosity changes with increasing shear rate, namely shear-thinning, where viscosity decreases, and shear-thickening, where viscosity increases.

Because of its importance in industry and in understanding nature, the study of Non-Newtonian Natural Convection (NNC) has attracted much research attention. Since closed-form solutions are rare and limited to extreme simplifications, the problem is usually treated numerically~\cite{yang2019comprehensive}. The dynamics of non-Newtonian fluids have been studied in various fields, from non-Newtonian blood flow in arteries~\cite{Kwack2014} to injection moulding of plastics \cite{Wang2010}. In particular, natural convection in cavities has been thoroughly analysed due to its direct application in industry~\cite{yang2019comprehensive}. 

From the numerical analysis point of view, the most commonly reported approaches to solving NNC in cavities are based on the Finite Volume Method (FVM) for the discretisation of the relevant partial differential operators and Semi-Implicit Method for Pressure-Linked Equations (SIMPLE) for pressure-velocity coupling~\cite{patankar1980numerical}. A comprehensive study of NNC solved with FVM and Quadratic Upstream Interpolation for Convective Kinematics (QUICK) to treat the convective terms was presented in~\cite{kim2003transient}, and with upwind stabilisation of convective terms in~\cite{turan2011laminar}. The FVM with QUICK was further investigated in~\cite{bozorg2019two} in solving NNC with internal rotating heater and cooler. Moraga et al.~\cite{moraga2017geometric} demonstrated the multigrid FVM with a fifth power differentiation scheme for convective transport in the solution of NNC with phase change. A similar numerical approach was also used in the solution of NNC in three dimensions (3D)~\cite{vasco2014parallel}. In conjunction with the Finite Difference Method (FDM) and upwind stabilisation, the solution of NNC with internal heat source was recently demonstrated in~\cite{loenko2019natural}. The NNC was also solved with a Finite Element Method (FEM)~\cite{alsabery2017transient, mishra2018natural} as well as with the Lattice Boltzmann Method (LBM)~\cite{kefayati2014simulation}.

All the aforementioned solution methods are mesh-based, i.e., the nodes are structured into polygons that completely cover the computational domain, a process also known as meshing. In FDM, FVM and LBM, a regular grid is often used, making meshing a trivial task, but at the cost of complications with irregular geometries and potential refinement. At FEM, meshing is mandatory and often also the most time-consuming part of the entire solution process, especially for realistic 3D geometries, which generally cannot be automatically meshed and therefore often require the user's help. An alternative to the mesh-based methods is the meshless approach~\cite{liu2002mesh}. The conceptual difference between mesh-based and meshless methods is that in the latter all relationships between nodes are defined solely by inter-nodal distances. An important implication of this distinction is that the meshless methods can work with scattered nodes, which greatly facilitates the consideration of complex 3D geometries and adaptivity. 

In this paper, we introduce a novel refined Radial Basis Function-generated Finite Difference (RBF-FD)~\cite{tolstykh2003using} meshless solution of NNC in two dimensional (2D) and 3D irregular domains computed on automatically generated scattered nodes~\cite{slak2019generation}, using the artificial compressibility method (ACM) for pressure velocity coupling~\cite{trojak2022acm, kosec2018localIrregularFlow, yasuda2023bulkViscosityACM}, and explicit time stepping. 
No stabilisation of the convective terms is used to minimise the impact of numerical diffusion in the results. We present a unified NNC solution procedure for 2D and 3D that can be easily extended to arbitrary geometries and inherently supports $hp$-adaptivity~\cite{slak2019adaptive, jancic2021p}.

Results for a reference 2D case are compared with data from Turan et al.~\cite{turan2011laminar} and Kim et al.~\cite{kim2003transient}, showing that our results are in the range of the comparative data. Therefore, in addition to a novel solution procedure, we also extend the range of available numerical solutions for the given problem with a completely different numerical approach. Note that the reference solutions differ only in the stabilisation of the convective terms (Upwind vs QUICK), otherwise in both papers results are computed with FVM and SIMPLE coupling.

In section~\ref{ch:problemFormulation} a mathematical model of NNC is discussed, followed by a presentation of the meshless numerical method and solution procedure in section~\ref{ch:numerical}. The analysis of results for the reference 2D case and a showcase of the method's versatility on more complex cases are shown in sections~\ref{ch:results} and \ref{ch:flexibility}.

\section{Problem formulation}
\label{ch:problemFormulation}

The dynamics of natural convection in non-Newtonian fluids are governed by a system of three partial differential equations describing the continuity of mass, the conservation of momentum, and the conservation of energy
\begin{align}
	\div \vec{v} &= 0, \label{eq:physics1}\\
	\rho (\pdv{\vec{v}}{t} + \vec{v} \cdot \grad{\vec{v}}) &= -\grad 
	p + 
	\div(\eta \left(\grad\vec{v} + (\grad{\vec{v}})^T\right)) -\vec{g} \rho \beta T_\Delta, \label{eq:physics2}\\
	\rho c_p (\pdv{T}{t} + \vec{v} \cdot \grad{T}) &= \div(\lambda 
	\grad
	T),\label{eq:physics3}\\
	\eta &= \eta_0 \left(\frac{1}{2}\norm{\grad{\vec{v}} + 
	(\grad{\vec{v}})^T}
	\right) ^{\frac{n-1}{2}},\label{eq:physics4}
\end{align}
with $\vec{v}$, $T$, $p$, $\eta$, $\rho$, $\vec{g}$, $\beta$, $T_\Delta$, 
$c_p$, 
$\eta_0$, $n$, representing the flow velocity field, temperature 
field, pressure field, viscosity field, density, gravity, thermal expansion coefficient, 
temperature 
offset, heat capacity, viscosity constant and non-Newtonian 
exponent, 
respectively.

The buoyancy force that drives the natural convection dynamics is 
relatively 
weak and ensures that the maximum velocity remains well below the 
speed of 
sound. This allows us to model the fluid as incompressible 
\cite{anderson2010fundamentals} and reduce the continuity equation 
to 
\eqref{eq:physics1}. The fluid motion is described by the 
Navier-Stokes 
equation \eqref{eq:physics2}, which is modified from its usual form 
by the 
addition of a force term describing the buoyancy caused by the 
thermal 
expansion. This force is approximated by the Boussinesq 
approximation 
\cite{Tritton1988}, which is based on the assumption that the 
acceleration of 
a fluid driven by natural convection remains insignificant compared 
to 
gravity and consequently the small thermal fluctuations of the 
density only 
play a role when amplified by the strong gravity in the buoyancy 
term. 
The Boussinesq approximation couples the fluid motion described by 
the 
Navier-Stokes equation \eqref{eq:physics2} with the temperature 
described by 
the energy equation \eqref{eq:physics3} and thus establishes 
the model 
for natural convection driven by the presence of a temperature 
gradient.

The constant viscosity $\eta$ is replaced by the Ostwald-de 
Waele 
power law model \cite{yang2019comprehensive} defined in Eq. 
\eqref{eq:physics4}. The shear dependence is captured by the 
exponentiated tensor 
norm of the shear rate tensor \footnote{The shear rate tensor for 
	incompressible fluids is the same as the strain rate tensor, 
	which can be 
	expressed as the symmetric part of the velocity gradient.} and 
	controlled by 
the exponent $n$, which controls the extent of the non-Newtonian 
behaviour, and 
$\eta_0$, which is used as a scaling factor.

\begin{figure}
	\includegraphics[width=\linewidth]{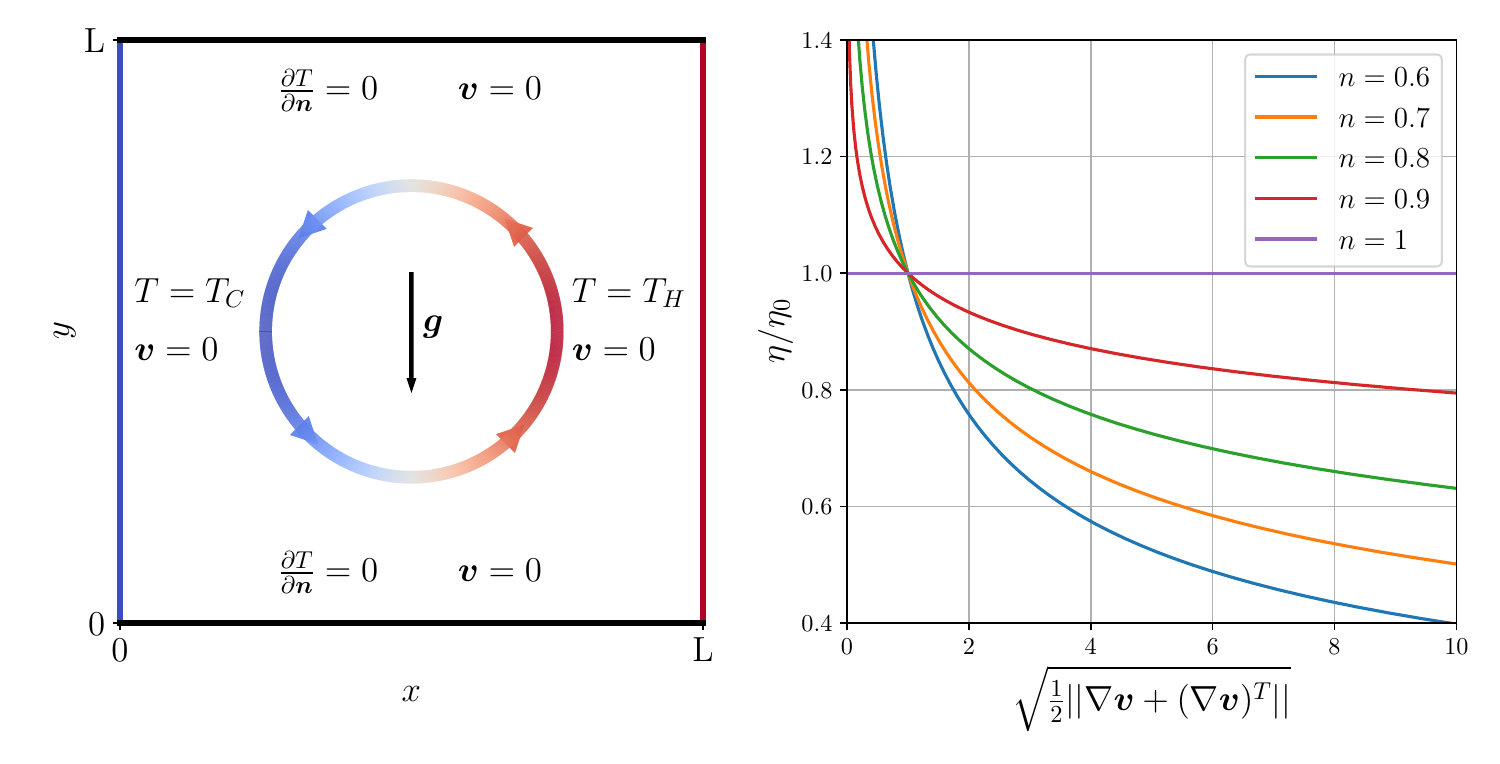}
	\caption{\textit{left:} A schematic representation of the De 
	Vahl Davis 
		differentially heated cavity case with velocity and 
		temperature boundary 
		conditions. \textit{right:} Relationship between viscosity 
		and shear rate 
		for a $n = 1$ Newtonian and a selection of shear-thinning $n 
		< 1$ 
		non-Newtonian cases.}
	\label{fig:caseDVD}
\end{figure}

The model reduces to a Newtonian fluid when $n=1$ and can be used to 
describe 
both shear thickening behaviour with $n > 1$ and shear thinning 
behaviour with 
$n < 1$. We focus on the latter as it leads to stronger convection, 
i.e. more 
interesting flow behaviour, and is more common in realistic fluids. 
A 
visualisation of the viscosity dependence for an arbitrarily chosen 
range of 
shear rate norms can be found on the right-hand panel of 
Figure~\ref{fig:caseDVD} for a range of exponents that we will use 
in further 
analysis. It can be seen from the figure that while $n$ has a 
dramatic impact 
on the effective viscosity, there is also a problem with the 
divergent 
viscosity at low shear rates in the shear thinning regime. This is 
mainly a 
problem for the initial, transient part of the simulation, where a 
flow pattern 
emerges from the initially stationary fluid and can be solved by 
downward 
bounding\footnote{The limit is arbitrarily chosen and has a 
negligible effect on 
	the result \cite{kim2003transient}.} the shear rate norm used in 
	the power law 
to $10^{-10}$. 

Model parameters that reflect the behaviour of a realistic fluid are typically  
determined by fitting the model to experimental data. Ostwald-de 
Waele power law model is the most basic approach to characterizing shear-dependent behaviour, and it is only appropriate for the intermediate shear rate regime. To accurately capture the asymptotic behaviour of realistic shear-thinning fluids at both low and high shear rates, the specialized models with more intricate algebraic structure due to additional parameters are used~\cite{Chhabra2010}. Nevertheless, from a numerical standpoint, these more complex models do not significantly affect the solution methodology nor do they introduce any additional numerical challenges.

The natural convection described with the system of equations 
described above 
is applied to the De Vahl Davis case \cite{DeVahlDavis1983}, a 
differentially 
heated square cavity with height and width $L = 1$, shown 
schematically in the 
left graph of Figure~\ref{fig:caseDVD}. The left wall is kept at a 
constant 
temperature $T_C = -1$, while the right wall is kept at a higher 
constant 
temperature $T_H = 1$, inducing the heat transfer that drives the 
dynamics of 
the system. The top and bottom boundaries are insulated. No-slip 
boundary 
conditions for velocity are imposed on all walls.

The difference in boundary temperatures creates a temperature 
gradient within 
the cavity that leads to natural convection. The convective flow 
appears due to 
the temperature-dependent changes in fluid density and the resulting 
buoyancy 
forces, which are described by the Boussinesq approximation term in 
Eq.~\eqref{eq:physics2}. The fluid cools and becomes denser at the 
left 
wall, so it falls and moves to the right, where it is heated by the 
hot wall, 
rises and completes the circular flow.

The circular flow caused by natural convection significantly 
increases the heat 
transfer between the differentially heated walls compared to 
conduction alone. 
The ratio between the two is known as the Nusselt number and 
provides a 
convenient reduction that expresses the convective heat transfer in a single 
scalar value. 
Such reduction can be further analysed to determine the temporal 
behaviour, as 
shown in Figure~\ref{fig:nusseltEvolution}, from which we can deduce 
when the 
steady state was reached. The Nusselt number used in the following 
analysis is 
calculated as the average of the values in the cold wall nodes
\begin{equation}
	\mathrm{Nu} = \frac{L}{T_H-T_C}\abs{\pdv{T}{x}}_{x=0}.
\end{equation}

The problem is further characterised by two dimensionless numbers. 
The Prandtl 
number (Pr) is a material property that expresses the ratio between 
the heat 
and momentum transport properties of the fluid. The Rayleigh number 
(Ra) 
can be interpreted as the ratio of buoyancy and thermal diffusivity, and is the product of the Grashof and Prandtl numbers. For the purposes of this study it can be interpreted as 
an analogue 
of the Reynolds number for natural convection, with larger values 
implying 
wilder dynamics. Both Ra and Pr are a function of viscosity and must 
be 
modified to account for non-Newtonian viscosity
\begin{align}
	\mathrm{Pr}=&\frac{\eta_0}{\rho} \alpha^{n - 2} L^{2 - 2 n}, \\
	\mathrm{Ra}=&\frac{\rho g \beta \Delta T L^{2 n + 1}}{\alpha^n 
		\eta_0},
\end{align}
with $\alpha = \frac{\lambda}{c_p \rho}$ as the thermal diffusivity. 
The definitions for the dimensionless numbers match the reference 
solution 
\cite{turan2011laminar} to facilitate comparison, as does the 
dimensionless time
\begin{equation}
	\hat{t} = \frac{\alpha}{L^2} t
\end{equation}
and velocity
\begin{equation}
	\hat{v} = \frac{L c_p \rho}{\lambda} v.
\end{equation}
We omit the Caret notation, since dimensionless values are used 
wherever 
velocity or time are referred to in the following discussion.

\section{Numerical solution procedure}
\label{ch:numerical}
Our goal is to solve the problem 
\eqref{eq:physics1}-\eqref{eq:physics3} using an 
approach that is as general as possible, including the generality of 
the number 
of dimensions, the order of the method, and the shape of the 
considered domain, 
while supporting a spatially variable discretisation. To achieve 
this, we 
employ the RBF-FD method with augmenting monomials. A meshless numerical 
technique 
operating on scattered nodes.

The first step is to discretise the domain, which in the meshless 
context means 
populating with scattered nodes. Although in the early stages of 
meshless 
development some authors used even randomly generated 
nodes~\cite{liu2002mesh}, 
it is now generally accepted that despite the apparent robustness of 
meshless 
methods regarding the node positioning, nodes still need to be 
generated 
according to certain rules~\cite{slak2019generation}, i.e. nodes 
have to 
''uniformly'' cover the domain with minimal empty spaces and satisfy 
the 
minimum distance requirement to avoid ill-conditioning of the 
approximation. 
There are several specially designed algorithms for meshless 
discretisation, 
ranging from expensive iterative~\cite{hardin2004discretizing} to 
advancing 
front~\cite{shankar2018robust} approaches. In this paper we use a 
Poisson disk 
sampling based advancing font approach~\cite{slak2019generation, 
	depolli2022parallel, duh2020fast, duh_discretization_2024}. The core of the algorithm is the iteration, where 
candidate nodes are sampled around the already positioned nodes and 
only 
those that do not violate the minimal distance requirement are added 
to the 
list of discretisation nodes. This conceptually simple approach has 
several 
convenient features. It is dimensionally agnostic, meaning that the 
formulation 
of the algorithm is the same regardless of the number of dimensions 
of the 
domain. It supports variable density node distributions, which are a 
cornerstone of the refined solutions discussed in 
\ref{ch:refinement}. It 
guarantees minimal spacing and is proven to be computationally 
efficient.

\begin{figure}
	\includegraphics[width=\linewidth]{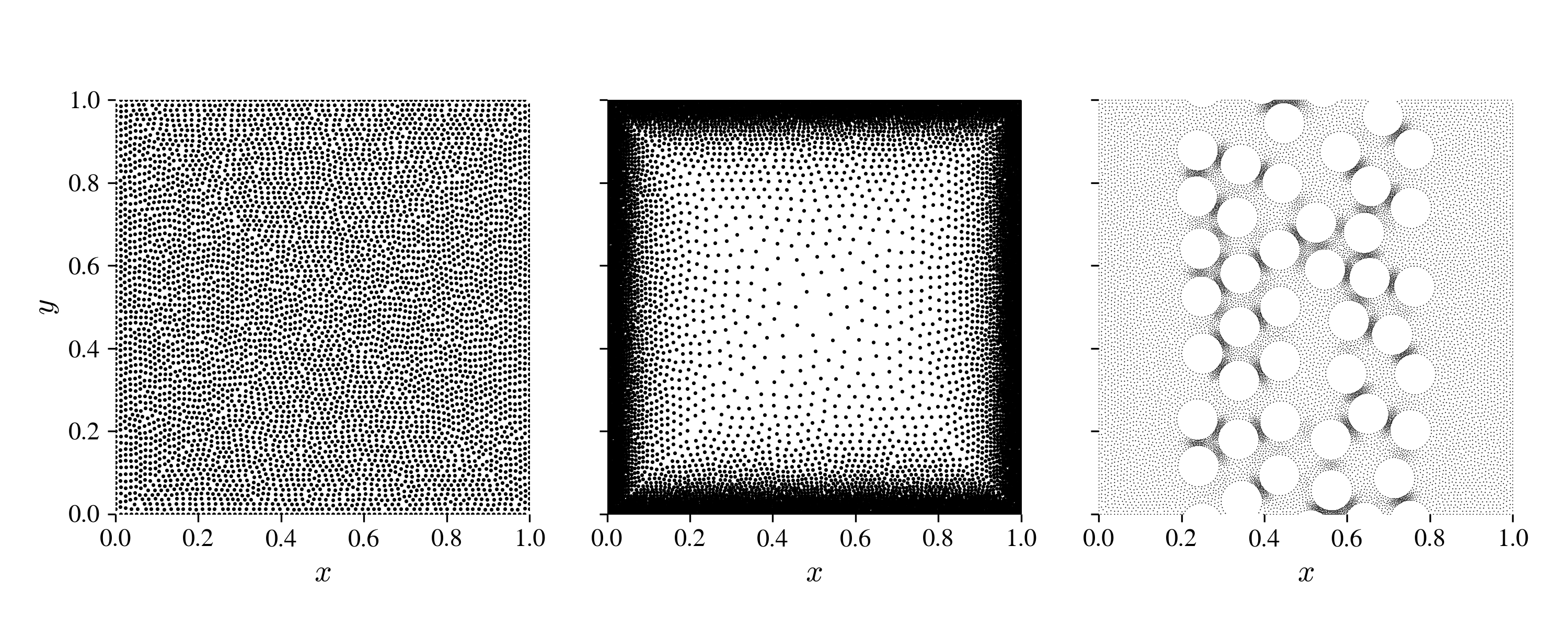}
	\caption{A comparison between node distributions in domains 
	populated with 
		a constant node density on the left, refined density in the 
		middle and the 
		obstructed domain with refined narrow channels on the right.}
	\label{fig:nodePlacement}
\end{figure}

We begin the analyses with a constant density node distribution, shown in the 
left panel of 
Figure~\ref{fig:nodePlacement}), which, as we will see in section 
\ref{ch:refinement}, is a sub-optimal strategy for the problem at hand 
with 
intense dynamics in the boundary layer (see section~\ref{ch:gfb}). An 
improved node placement strategy, shown in the central panel of 
Figure~\ref{fig:nodePlacement} and discussed in section~\ref{ch:refinement}, 
uses prior 
knowledge of the system to devise a target node density that 
significantly 
reduces computational cost without sacrificing accuracy. This 
solution is still 
not ideal, as it requires physical intuition and manual input, but 
it 
demonstrates the potential for the eventual goal, which is an 
$h$-adaptive 
solution with an appropriate error indicator. In a rightmost plot of 
Figure~\ref{fig:nodePlacement} we demonstrate discretisation of a an 
irregular 
domain that will be used in section~\ref{ch:flexibility}.

Once the computational nodes are generated we identify the approximation stencil, i.e. the nodes used in local approximation of differential operators. Here, we will resort to simplest stencil strategy, where for each node $i$, $s$ neighbouring nodes are identified as a stencil $S_i$. Besides basic closest-node stencils, recent studies \cite{davydov2023stencilSelection} suggest that advanced symmetric stencils might significantly enhance the stability of meshless approximations. However, most research still employs nearest neighbours stencils; thus, we will also use it in this paper. 


In the next step we construct a generalised finite difference 
approximation to 
numerically evaluate the linear differential operator $\L$ in the 
central node 
from the approximated function's values in stencil nodes
\begin{equation}
	(\L u)_i \approx \sum_{j=1}^{s} w_{i, j} u_{S_i(j)},
	\label{eq:operatorApprox}
\end{equation}
where $u_k$ denotes the value of the arbitrary approximated function 
$u$ at 
the position of the $k$-th node $\vec{p}_k$. The weights $w$ are 
determined through the demand that the Eq. \eqref{eq:operatorApprox} 
is 
exact for a set of basis functions, in our case radial basis 
functions 
(RBF) 
\begin{equation}
	\label{eq:rbfScaling}
	\phi(j, k) = \phi \left( \frac{\norm{\vec{p}_k - 
	\vec{p}_j}}{\delta_j} 
	\right),
\end{equation}
where we introduced $\delta_j$ as a local scaling factor that 
decouples the 
approximation from the choice of coordinate system. It  is set to an 
arbitrary 
local measure of distance, e.g., the distance to the closest stencil 
node, and 
is of the utmost importance when using RBFs that use a scaling 
parameter.  
Finally, we get a local linear system $\vec{A} \vec{w}_i = \vec{b}$
\begin{gather}
	\label{eq:nonAugmentedSystem}
	\begin{split}
		\begin{bmatrix}
			\phi(S_i(1), S_i(1)) & \cdots  & \phi(S_i(1), S_i(s)) \\
			\vdots & \ddots & \vdots \\
			\phi(S_i(s), S_i(1)) & \cdots & \phi(S_i(s), S_i(s)) \\
		\end{bmatrix}
		\begin{bmatrix}
			w_{i, 1} \\ \vdots \\ w_{i, s}
		\end{bmatrix}
		=
		\begin{bmatrix}
			(\L\phi)(i, S_i(1)) \\
			\vdots \\
			(\L\phi)(i, S_i(s))
		\end{bmatrix},
	\end{split}
\end{gather}
for each node with the solution yielding stencil weights 
$\vec{w}_i$. The right 
hand side vector $\vec{b}$ is formed by applying the linear operator 
$\L$ to 
the basis function and evaluating the result with an argument 
analogous to 
Eq.~\eqref{eq:rbfScaling}. 
We use polyharmonic 
splines (PHS) basis
\begin{equation}
	\phi(r) = r^{k},
\end{equation}
with odd order $k$, additionally augmented with polynomials for ensuring the polynomial reproduction and positive definiteness~\cite{flyer2016polynomials1}. This setup is widely recognised as an RBF-FD method in meshless community. 


 
The system for stencil weights is expanded with $N_p = \binom{m + 
d}{m}$ 
monomials\footnote{The 6 monomials in 2-D case with $m=2$ would be 
$q = \{1, 
	x, y, x^2, xy, y^2\}$.} $q$, where $m$ denotes the monomial 
	order and $d$ the 
spatial dimension. The monomials are scaled with a similar 
argumentation as 
RBFs in Eq.~\eqref{eq:rbfScaling}
\begin{equation}
	q_l(j, k) = q_l \left( \frac{\vec{p}_k - \vec{p}_j}{\delta_j} 
	\right).
\end{equation}
The $\vec{A} \vec{w}_i = \vec{b}$ system for approximation weights 
from 
Eq.~\eqref{eq:nonAugmentedSystem} is augmented with monomials
\begin{equation}
	\begin{gathered}
		\begin{bmatrix}
			\vec{A} & \vec{Q} \\
			\vec{Q}^T & 0 \\
		\end{bmatrix}
		\begin{bmatrix}
			\vec{w}_i \\ \vec{\lambda}
		\end{bmatrix}
		=
		\begin{bmatrix}
			\vec{b} \\
			\vec{c} \\
		\end{bmatrix},
		\\
		\vec{Q} = \begin{bmatrix}
			q_1(S_i(1), S_i(1)) & \cdots  & q_{N_p}(S_i(1), S_i(1))\\
			\vdots & \ddots & \vdots \\
			q_1(S_i(s), S_i(1)) & \cdots & q_{N_p}(S_i(s), S_i(1))\\
		\end{bmatrix}
		,
		\vec{c} = \begin{bmatrix}
			(\L q_1)(S_i(1), S_i(1)) \\
			\vdots \\
			(\L q_{N_p})(S_i(1), S_i(1))
		\end{bmatrix},
	\end{gathered}
\end{equation}
with the additional weights $\vec{\lambda}$ treated as Lagrange 
multipliers 
and discarded after computation.

Augmentation with an order of at least $m = \frac{k - 1}{2}$ is 
required to 
guarantee the positive definiteness for a PHS with order $k$. Higher 
orders of augmentation $m$ 
provide better convergence characteristics \cite{jancic2021monomial, leborne2023guidelinesRBFFD}, with the order of convergence
\begin{equation}
	\mathcal{O}(h^{m + 1 - \ell})
	\label{eq:expectedOrder}
\end{equation} 
for $m \geq l$, where $\ell$ is the order of the linear operator $\L$. The higher accuracy comes at the 
cost of increased computational complexity, since the required stencil size is 
$s >= N_p$, with $s > 2N_p$ as the often recommended value 
\cite{bayona2017}. Increased stencil size affects both 
the 
pre-computation of the approximation weights, with complexity 
$\mathcal{O}((s + 
N_p)^3)$, and the eventual scalar product $\mathcal{O}(s)$ 
evaluation. We will 
use $k=3$, which is the minimum odd value that can be used to approximate second order derivatives, for all numerical results presented in this paper and compare the results for monomial orders $m=2$ and $m = 4$. The support size is chosen conservatively as $s = 2N_p + 1$ unless otherwise specified.  Nevertheless, augmenting monomials open an opportunity for $hp$-adaptivity, since their order directly controls the order of the method~\cite{jancic_strong_2023}.
 
Now that we have the RBF-FD approximation for derivatives, next step is to formulate a solution procedure for the problem at hand. The pressure-velocity coupling in the incompressible Navier-Stokes system is performed using the artificial compressibility method (ACM)\cite{trojak2022acm, kosec2018localIrregularFlow, yasuda2023bulkViscosityACM} that was first introduced by Chorin~\cite{chorin1967acm} in 1967 but is now experiencing resurgence due to its explicit and local nature that allows for easy parallelisation and GPU usage~\cite{kosec_super_2014,kajzer2018GPU}. The method works by artificially introducing a slight compressibility into the otherwise incompressible system in order to calculate the pressure field. 

Since the focus of this study lies in the spatial discretisation we use the simple first-order explicit Euler method for the temporal discretisation of the Navier-Stokes equation. The intermediate velocity 
\begin{equation}
	\vec{v'} = \vec{v} + \Delta t \left( \div(\eta \left(\grad\vec{v} + (\grad{\vec{v}})^T\right)) - \vec{v} \cdot \nabla\vec{v} -\vec{g} \rho \beta T_\Delta \right) \label{eq:MNS_intermediate}
\end{equation}
is calculated by using the viscosity field $\eta$ calculated from the previous step velocity according to \ref{eq:physics4}, and the offset $T_\Delta$ in Bousinessq term based on the previous step temperature field. Note that the pressure term is omitted at this step.
The time-step
\begin{equation}
	\Delta t = \min\left(\min_i\left(c_1\frac{h_i}{\lVert\vec{v_i}\rVert_2}\right), \min_i\left(c_2\frac{\rho h_i^2}{2 \eta_i}\right)\right)
\end{equation}
is calculated dynamically during each step of the iteration to adapt to the changing velocity and viscosity in computational nodes with a wide range of inter-nodal distances $h_i$. The time-step constants are chosen as $c_1 = 0.05$ and $c_2=0.15$ based on trial and error. Values of the dynamic timestep during the simulation with a constant discretisation density using $h=0.005$ are shown on the right graphs of 
Figure~\ref{fig:nusseltEvolution}. The variable time-step is beneficial for most of the displayed cases as the most stringent requirements only occur during the relatively short period of initial flow formation. Note that the $h^2$ term in the viscous limit is likely to dominate the timestep selection for very dense discretisations and might prove to be a limiting factor.

We omitted the pressure term from the initial intermediate velocity $\vec{v'}$ because the subsequent pressure-velocity coupling is done iteratively. First the intermediate velocity is corrected with the previous pressure field
\begin{equation}
	\vec{v} = \vec{v'} - \frac{\Delta t}{\rho} \nabla p, \label{eq:MNS_velocity},
\end{equation}
then the pressure field is corrected to counteract any divergence present in the velocity
\begin{equation}
	p  \leftarrow p - \Delta t C^2  \rho  (\nabla \cdot \vec{v}), \label{eq:MNS_pressure}
\end{equation}
while enforcing
\begin{equation}
	\pdv{p}{\vec{\hat{n}}} = 0,
\end{equation}
on the boundary, with $\vec{\hat{n}}$ being the corresponding surface normal vector.
 The strength of the artificial compressibity effects is recomputed at every step and governed with the artificial speed of sound $C$~\cite{rahman2008ACM}
\begin{equation}
	C = \gamma \max(\max_i(\lVert\vec{v_i}\rVert_2), \lVert\vec{v}_{ref}\rVert_2), \label{eq:MNS_C}
\end{equation}
where $\gamma$ is the compressibility parameter, and $\vec{v}_{ref}$ a reference velocity introduced to prevent potential issues caused by $C$ reaching zero. We use $\gamma = 5$ for all computations presented in this paper.
If a time-accurate solution was required, the pressure-velocity coupling iteration would have to be repeated until the maximum divergence of the velocity field dropped below a desired level. In steady state cases like the DVD in examined regime a single iteration is sufficient as we are not really interested in the transient regime.

The new velocity is then used in the advection part of the energy transfer equation, again discretised with the explicit Euler method, which is then used to calculate the new temperature field. The same dynamic timestep $\Delta t$ calculated for the momentum equation can also be used for the energy equation as the latter is well within stability limits of the former when dealing with $\mathrm{Pr} =100 \gg 1$ case. 


\section{Results}
\label{ch:results}
\subsection{General flow behaviour}
\label{ch:gfb}
We will study the impact of non-Newtonian behaviour on three flow cases specified with Rayleigh and Prandtl dimensionless numbers as defined in section~\ref{ch:problemFormulation}. Rayleigh numbers $\mathrm{Ra} =\{10^4, 10^5, 10^6\}$ are chosen to capture changes in the upper range of steady natural convection regime, while $\mathrm{Pr}=100$ remains constant and is chosen to facilitate comparison with the existing reference solution \cite{turan2011laminar}. The non-Newtonian behaviour is examined with five values for the non-Newtonian 
exponent $n = \{0.6, 0.7, 0.8, 0.9, 1\}$ progressing from the $n=0.6$ case that exhibits the strongest shear-thinning behaviour to the 
Newtonian $n=1$ case. Note that the most extreme case with $\mathrm{Ra} = 10^6$ and $n = 0.6$ is close to the edge of the steady regime. Further increases in Ra or decreases in $n$ would result in an oscillatory flow~\cite{kim2003transient}.

\begin{figure}
	\includegraphics[width=\linewidth]{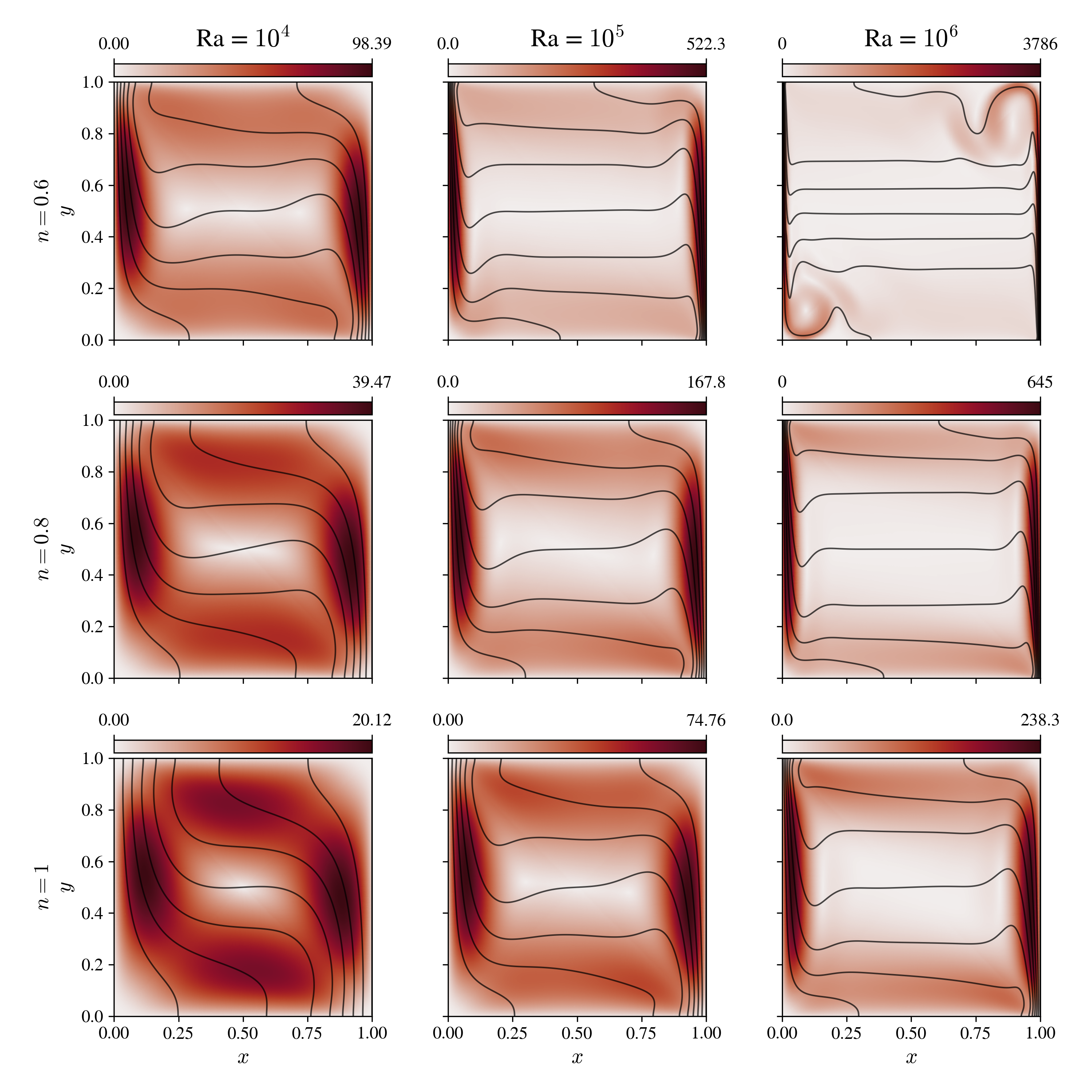}
	\caption{Flow profiles for a selection of cases. Velocity 
	magnitude is 
		visualised with a heat-map 	while the overlaid contours 
		display the 
		changes in 	temperature. Each sub-figure has a distinct 
		velocity range 
		specified by the colourbar above.
		}
	\label{fig:flowFields}
\end{figure}

A sample of the resulting flow profiles can be seen in 
Figure~\ref{fig:flowFields} where the velocity magnitude is 
displayed as a heat 
map and overlaid with temperature contours. All cases exhibit the 
previously 
described circulation caused by natural convection, but there are 
drastic 
differences in the maximum velocity and the thickness of boundary 
layers as Ra 
and shear-thinning increase.

Both effects are expected and can be explained by the definition of 
varied 
parameters. The Rayleigh dimensionless number is defined as the 
product of 
Prandtl, which we keep constant, and Grashof dimensionless numbers. 
The latter 
expresses a ratio between buoyant and viscous forces and explains 
why 
increased Ra results in cases where fluid in the boundary layer 
convects away 
before conducting much heat to the neighbouring fluid. Similarly 
reduced $n$ 
decreases viscous penalty for high velocity gradients close to the 
constant 
temperature edges, leading to a further reduction in boundary layer 
thickness.

\begin{figure}
	\includegraphics[width=\linewidth]{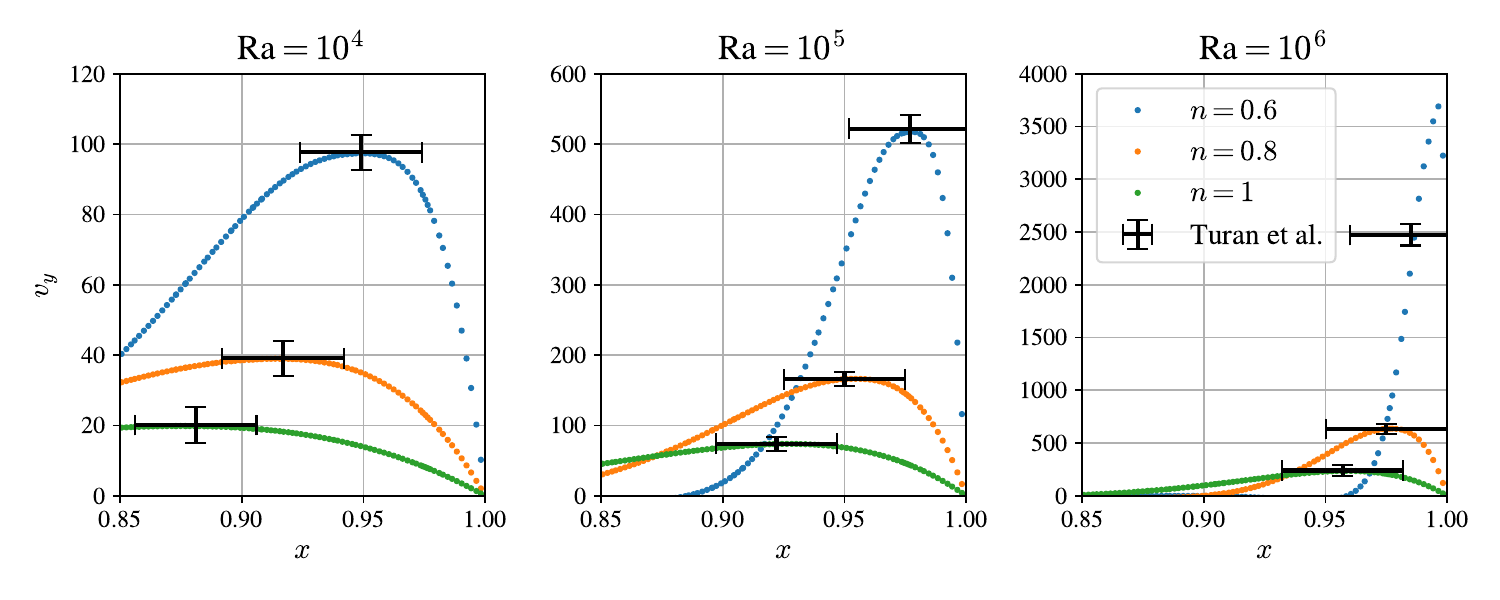}
	\caption{Cross-section of vertical velocity $v_y$ close to the right wall at $y = 0.5$ for a range of different Ra and $n$, calculated with $\mathrm{Pr}=100$ and $m=4$. The error bars show the locations and values of maximum vertical velocity in corresponding plots from a reference solution~\cite{turan2011laminar}.}
	\label{fig:velocityProfileRefComparison}
\end{figure}

Velocity fields also offer the first opportunity to verify the results. We compare the vertical velocity cross-section at $y=0.5$ with the reference solution~\cite{turan2011laminar} in Figure~\ref{fig:velocityProfileRefComparison}. The position and value of the largest vertical velocity in the reference solution is added with an estimated error caused by the process of extracting the data from a similar cross-section plot. Our results match the reference solution across the whole range of the considered parametres with the exception of the $\mathrm{Ra}=10^6$, $n=0.6$ case with problematic convergence that is discussed and solved in the following sections.

\begin{figure}
	\includegraphics[width=\linewidth]{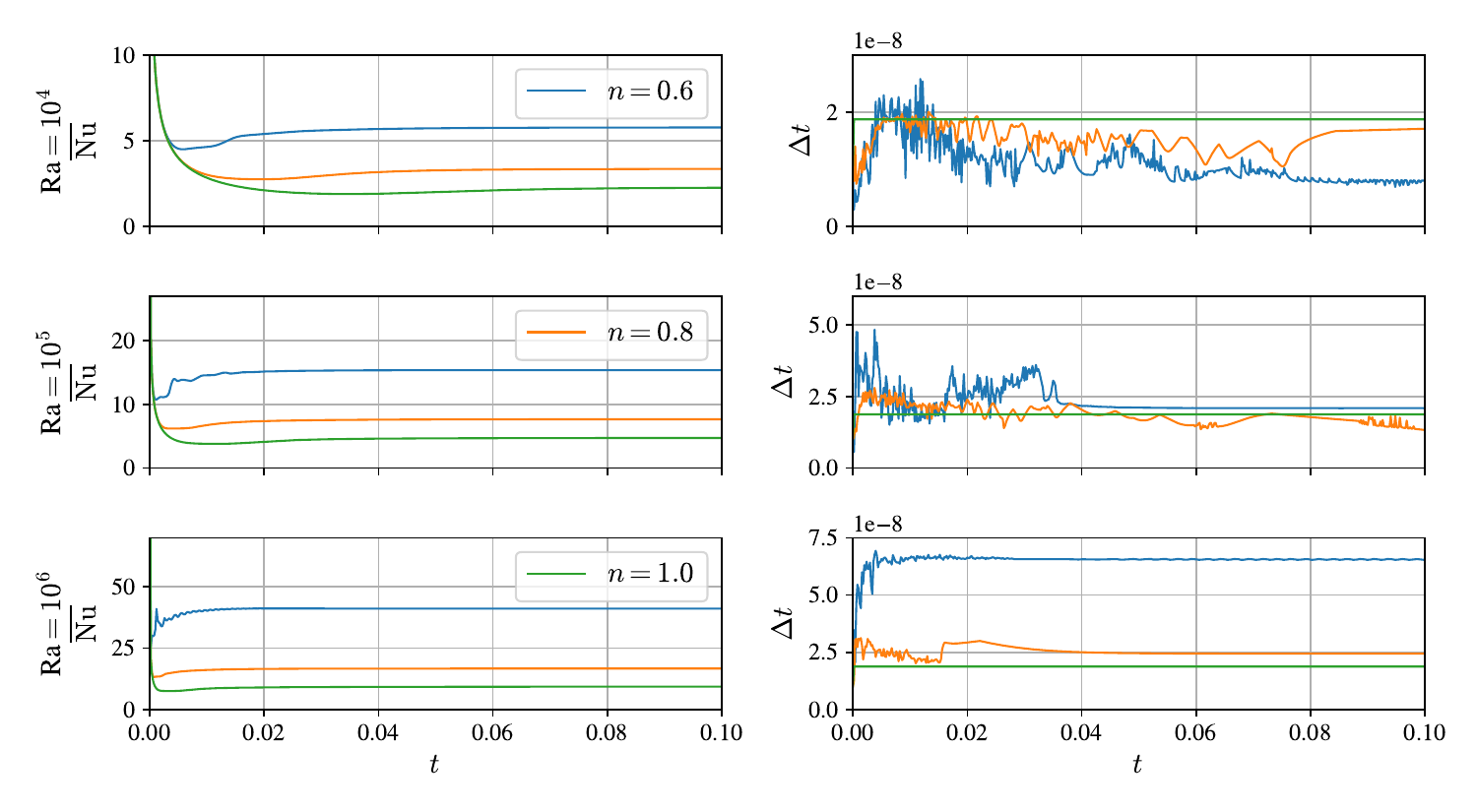}
	\caption{\textit{left:} Time evolution of the average Nusselt number on the cold wall that is used as a scalar observable for the system's dynamics. \textit{right:} Dimensionless values of the dynamic timestep throughout the simulation. Cases shown in this Figure use a constant discretisation density with internodal distance $h = 0.005$ corresponding to $N=35222$ computational nodes.}
	\label{fig:nusseltEvolution}
\end{figure}

In the subsequent analysis we utilise the Nusselt number as a scalar observable to simplify the description of the system's behaviour, facilitating temporal observation and comparison between different cases. One such example is shown in  
Figure~\ref{fig:nusseltEvolution} where we track the evolution of 
the Nusselt 
number to determine when we reach a steady state. We can also 
observe the 
effects of stronger shear-thinning. Cases with lower $n$ are faster 
to reach 
the stationary state as it is easier for convection to start, with 
the lower 
effective viscosity, and to play a bigger part in heat transfer as 
reflected in 
higher Nusselt values. The initial 
Nusselt number is 
high due to the high temperature gradient at the constant 
temperature boundary 
when starting from a zero-temperature zero-velocity initial 
condition. The 
value then decreases as the importance of conductive heat transport 
increases 
until circulation is established. Alternatively we could start with the diffusive field as the initial condition for temperature but this would only lead to a more violent flow formation with little benefit to accelerating the convergence to the convection dominated steady state. 

\subsection{Convergence under $h$ refinement}
The node density convergence analysis is performed to analyse node independence. We repeatedly run the same cases with decreasing internodal distance $h$ and examine how the average Nusselt number on the cold boundary changes with increasing number of computational nodes. The velocity and temperature fields are initially set to zero and the simulation is ran to $t=0.1$ which is determined to be sufficient to ensure that a stationary state has been reached as seen in Figure~\ref{fig:nusseltEvolution}.

\begin{figure}
	\includegraphics[width=\linewidth]{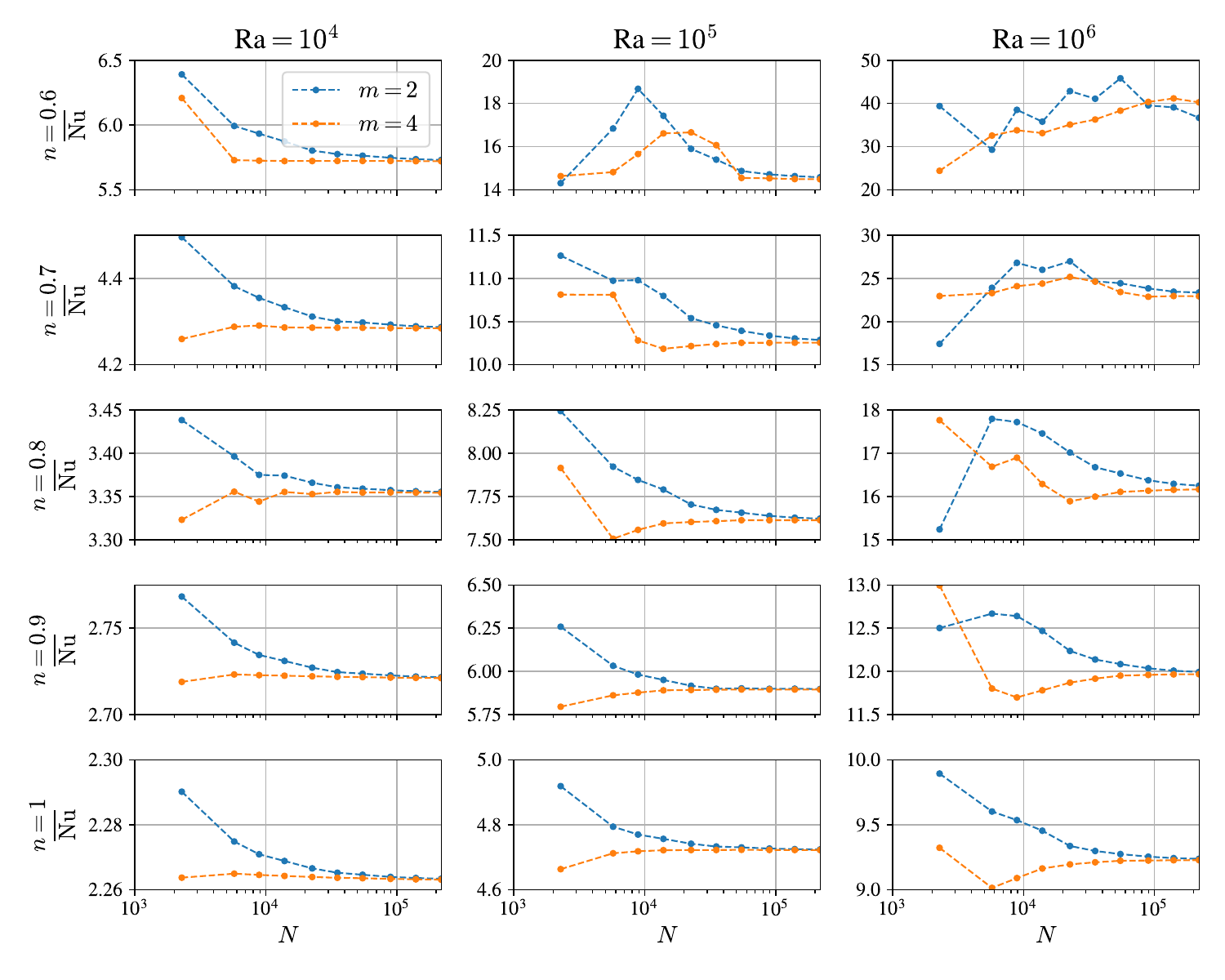}
	\caption{Convergence of the average Nusselt number on the cold 
		boundary 
		with columns for all of the considered Ra and rows for 
		different 
		non-Newtonian exponents $n$. The two different considered monomial augmentation orders 
		$m$ are shown 
		with different plot colours. Node count is proportional to 
		the inter-nodal 
		distance as the density is constant throughout the domain.}
	\label{fig:nusseltConvergence}
\end{figure}

The convergence behaviour for Newtonian ($n = 1$) cases shown in the bottom 
row of Figure~\ref{fig:nusseltConvergence} is very tame, as the coarsest considered
discretisation already provides results with less than 10\% discrepancy compared to the 
finest. 
Unfortunately, this is no longer the case as we progress towards 
non-Newtonian 
cases with a thinner boundary layer. The variation of the observed 
values 
increases and the convergence rate decreases. This culminates in the 
most 
extreme case with $\mathrm{Ra}=10^6$ and $n=0.6$, shown in the upper right 
graph of 
Figure~\ref{fig:nusseltConvergence}, where adequate convergence is 
not achieved with the considered node counts. This can be further corroborated with the reference mismatch already observed for this case in Figure~\ref{fig:velocityProfileRefComparison}. Proceeding to even higher node densities, 
requiring 
longer computational times, is wasteful, especially as the utilised 
numerical 
method allows for an elegant optimisation described in following 
sections.

\begin{figure}
	\includegraphics[width=\linewidth]{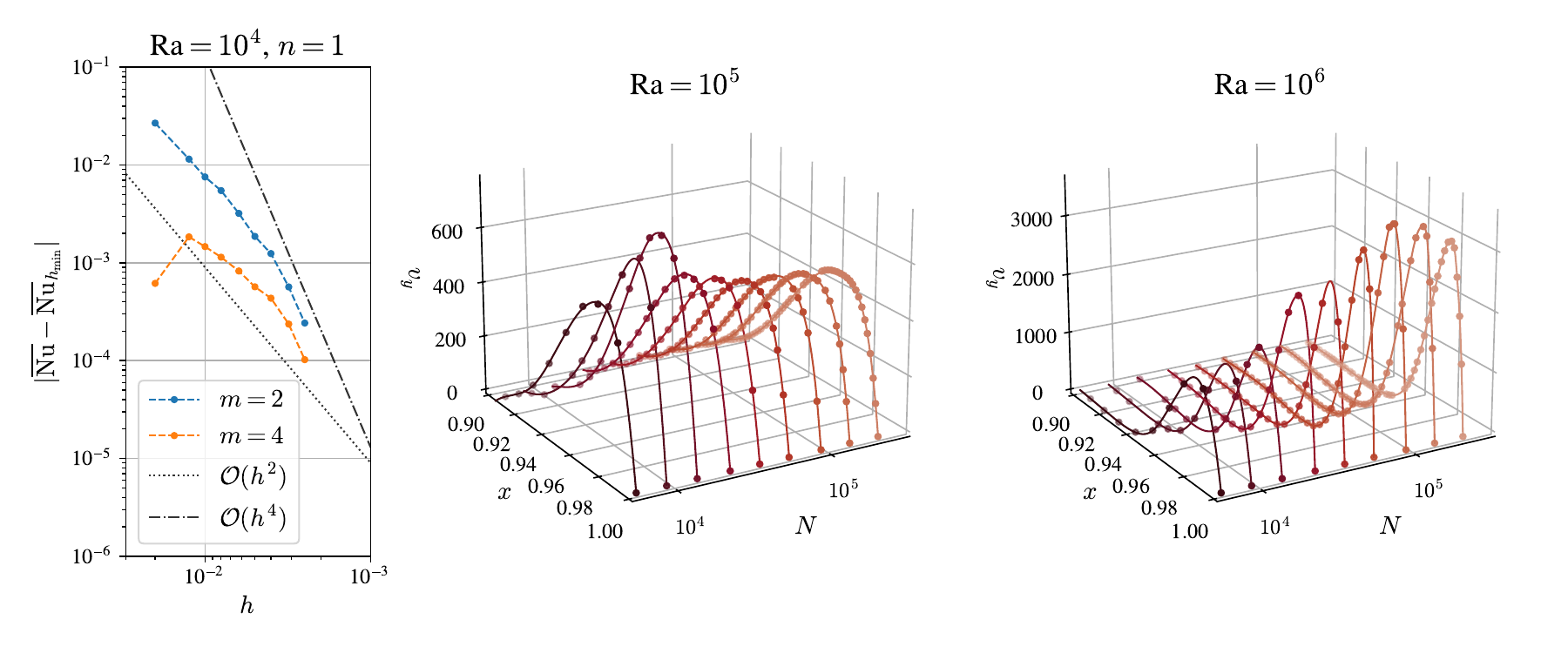}
	\caption{\textit{left:} Convergence of average Nusselt number on the cold boundary towards the value calculated with the finest discretisation that was used as the ground truth. \textit{middle \& right:} Changes in the interpolated vertical velocity profiles close to the boundary at $y = 0.6$ with increasing number of uniformly positioned computational nodes. Cases shown are parametrised with $\mathrm{Pr}=100$ and $n=0.6$ and computed with $m=2$. Dots display the horizontal position and vertical velocity in computational nodes that lie less than $\frac{h}{2}$ from the interpolation $y$.}
	\label{fig:logConvergenceProfiles}
\end{figure}

In next step we compare different monomial augmentation orders $m=2$ and $m=4$, shown as different colours in Figure~\ref{fig:nusseltConvergence}. 
We consider the finest available discretisation to produce an ''accurate`` solution and study the rate of convergence relative to that in the left graph of Figure~\ref{fig:logConvergenceProfiles}. From the log-log plot of the average Nusselt number offset against the inter-nodal distance $h$, we can determine that even though the higher augmentation order leads to a smaller initial error, the rate of convergence is similar for both. The discrepancy from the expected order, expressed with Eq.~\eqref{eq:expectedOrder}, are most likely caused by non-linearities in the system.

There is an interesting transition in convergence curves shown in 
Figure~\ref{fig:nusseltConvergence} as we move towards the upper right 
corner, 
corresponding with cases that exhibit a thinner boundary layer with 
higher 
velocities. The calculated Nusselt numbers initially rise as we 
increase 
node density before falling towards the value they eventually converge 
to. This 
behaviour can be better understood by examining the velocity 
cross-section 
convergence, shown in the middle and right graph of Figure~\ref{fig:logConvergenceProfiles}, 
for two of the $n=0.6$ cases where it is most apparent. The velocity profiles show that the boundary layer initially becomes narrower and faster until a sufficient number of nodes is present to adequately capture the dynamics. This behaviour provides us with additional motivation for refining the discretisation close to the heated and cooled boundaries. Furthermore the hypothesized issues with derivative approximation appear to also have a connection with the stencil size and/or approximation order as the Nusselt number peaks visible in Figure~\ref{fig:nusseltConvergence} move towards higher densities when augmented with a higher monomial order. 

\begin{figure}
	\includegraphics[width=\linewidth]{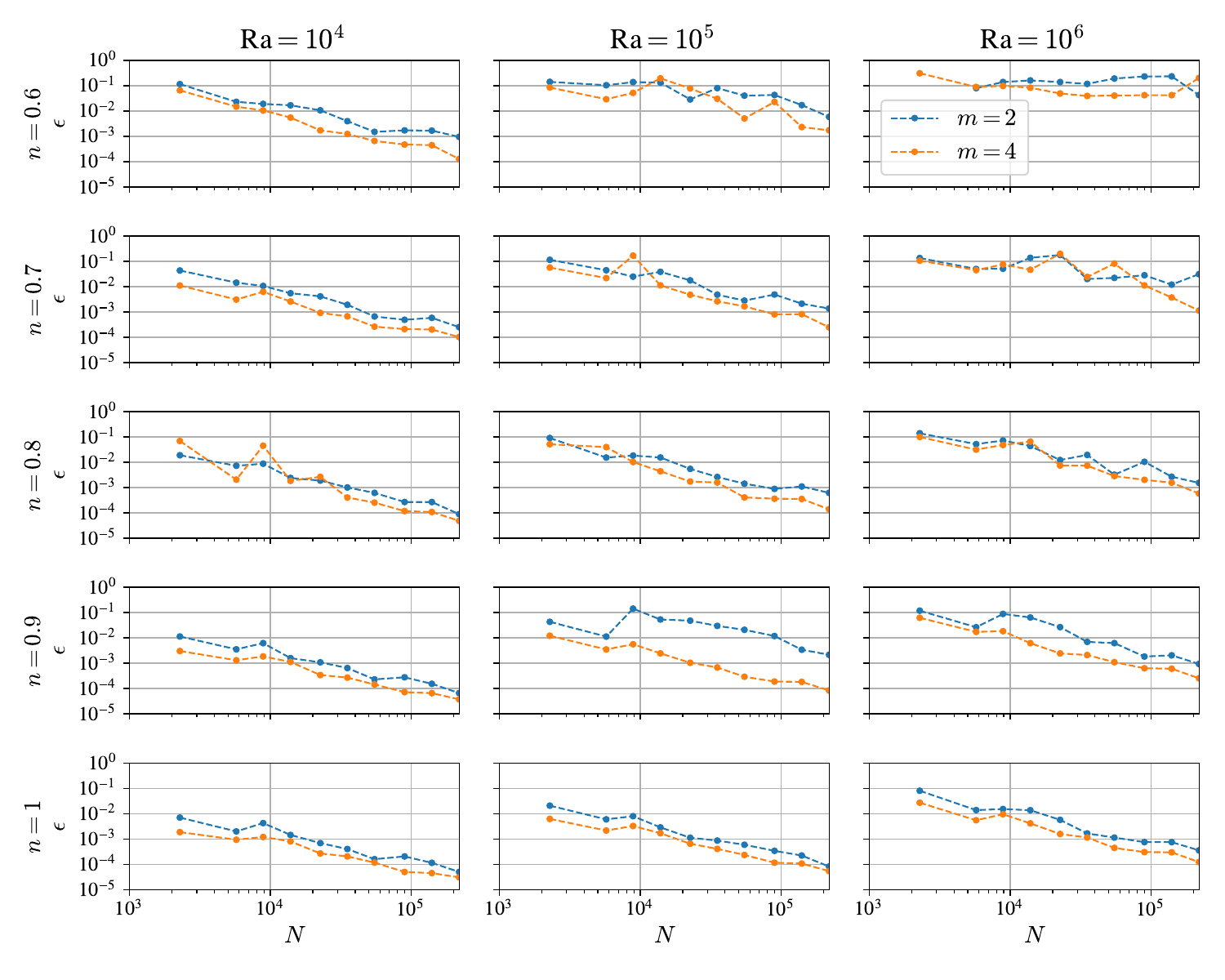}
	\caption{Convergence of the relative symmetry error at $y=0.75$ 
	with 
		columns for all of the considered Ra and rows for different 
		non-Newtonian 
		exponents $n$. The two different considered monomial augmentation orders 
		$m$ are shown 
		with different plot colours. Node count is proportional to the inter-nodal 
		distance as the 
		density is constant throughout the domain.}
	\label{fig:symmetryErrorGrid}
\end{figure}

Furthermore, we examine the violation of symmetry
\begin{equation}
	u(x, y) = -u(1-x, 1-y)
\end{equation}
as an alternative method for assessing the fitness of the solution. This is additionally motivated by the fact that the $\mathrm{Ra}=10^6$, $n=0.6$ case with problematic convergence exhibits clearly visible asymmetry in velocity and temperature profiles shown in Figure~\ref{fig:flowFields}.  
First, we 
introduce a measure of symmetry violation by interpolating the 
vertical 
velocity 
$v_y$ at an arbitrarily chosen $y$ and $1 - y$ and calculating the 
relative 
error as
\begin{equation}
	\epsilon = \frac{\max_{x \in [0, 1]}(\abs{v_y(x, y) + v_y(1-x, 
			1-y)})}{\max_{x \in [0, 1]}(\abs{v_y(x, y)})}.
\end{equation}
where we normalise the maximum offset between a vertical velocity 
and its 
symmetric value with the maximum vertical velocity for a given case 
and 
selected $y$. This choice of denominator is preferable as it focuses 
on the 
symmetry errors in the relevant high velocity part of the domain 
while 
remaining relative for comparison between different cases. The 
resulting 
symmetry errors as a function of the number of nodes are shown in 
Figure~\ref{fig:symmetryErrorGrid} and corroborate the previous 
discussion that 
was based on convergences in Figure~\ref{fig:nusseltConvergence}. Additionally, we can also confirm the observation from the left graph of Figure~\ref{fig:logConvergenceProfiles} as both $m=2$ and $m=4$ exhibit similar convergence behaviour.
Note that in 
addition to confirming the convergent behaviour of the method, the 
introduced 
symmetry violation can also be interpreted as a lower bound for the 
error.

\subsection{Refinement}
\label{ch:refinement}
Small flow structures and large velocity gradients (boundary layer) only occur in a 
relatively small part of the domain as seen in the right graph of Figure~\ref{fig:logConvergenceProfiles} whose $x$-axis only covers about 10\% of the domain's width. 
Needlessly covering the entire 
domain with the 
high node density, that is only required close to the cold and 
hot 
boundary, drastically increases the number of operations required at 
every 
time-step without improving the results.

\begin{figure}
	\includegraphics[width=\linewidth]{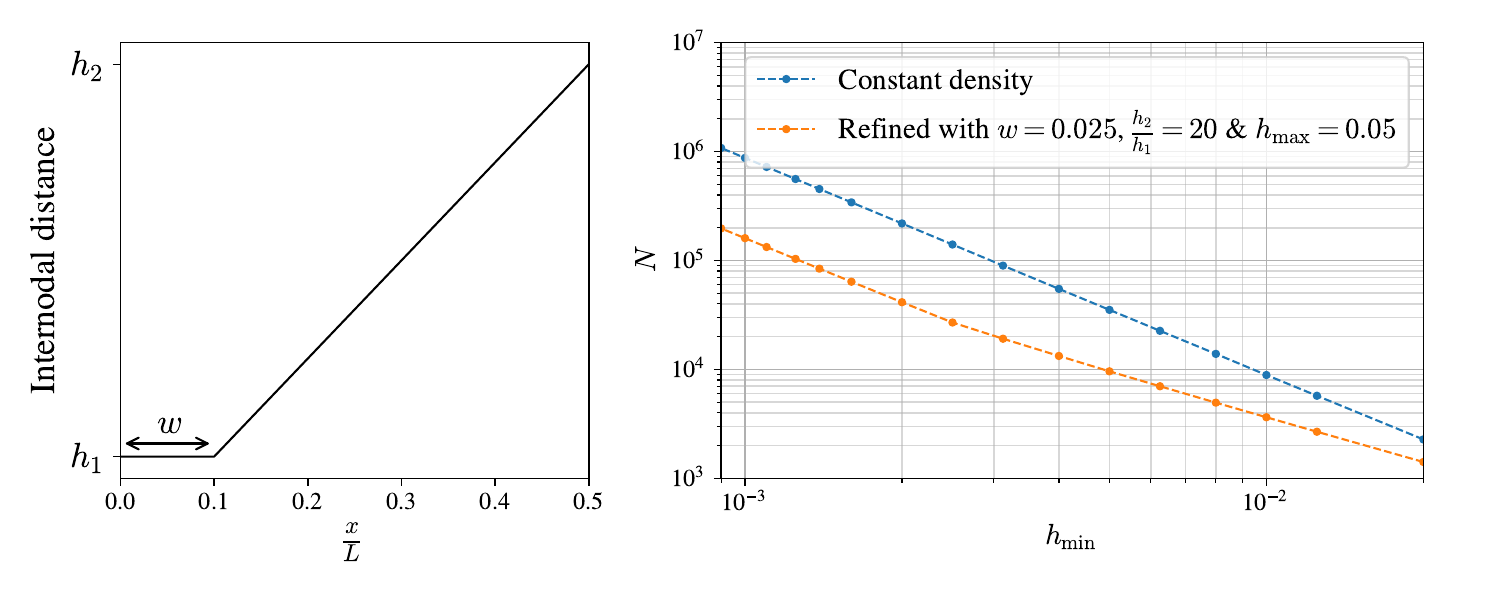}
	\caption{\textit{left:} A schematic representation of the inter-nodal 
	distance as 
		a function of distance from the boundary. \textit{right:} 
		Computational 
		node count as a function of minimum inter-nodal distance for 
		a constant 		
		density and a refined case.}
	\label{fig:refinementScheme}
	
\end{figure}

We introduce a variable node density expressed as the inter-nodal distance. In general, the inter-nodal distance would be a 
function of 
position but we use a symmetric configuration, shown schematically 
in the left 
graph of Figure~\ref{fig:refinementScheme}, that only depends on the 
distance 
to the closest boundary $d$. The idea is to position nodes with 
inter-nodal 
distance $h_1$ in a band within $w$ of the boundary to ensure a 
sufficient 
discretisation in this intense region, while using much coarser 
discretisation $h_2$ in the centre of the domain, with a linear 
transition between the 
two
\begin{equation}
	h=\left\{\begin{array}{ll}
		h_1, & \text{where } d < w \\
		h_1 + \frac{d - w}{\frac{L}{2} - w} (h_2 - h_1). & \text 
		{otherwise}
	\end{array}\right.
\end{equation}
The value of the central inter-nodal distance
\begin{equation}
	h_2 = \min(k_{\mathrm{ref}} h_1, h_\mathrm{max}),
\end{equation}
where $k_{\mathrm{ref}}$ is the refinement ratio that will be expressed as $\frac{h_2}{h_1}$ in subsequent analysis, is expressed in terms of minimal distance $h_1$ to prevent excessive density gradients that would require further analysis in terms of method stability and accuracy. Furthermore $h_2$ is bounded with $h_\mathrm{max}$ to prevent instabilities.
This relatively simple refinement scheme leads to drastic savings in node 
count shown in the right graph of Figure~\ref{fig:refinementScheme}. Even with the 
relatively conservative refinement parameters, there are almost an order of 
magnitude fewer nodes for the same boundary inter-nodal distances $h_1$. The slight non-linearity in the node count is caused by $h_\mathrm{max}$ limited $h_2$ when $h_1$ is large. The savings 
appear to be even more dramatic when a refined discretisation is visually 
compared to an unrefined one as seen in Figure~\ref{fig:nodePlacement}.

\begin{figure}
	\includegraphics[width=\linewidth]{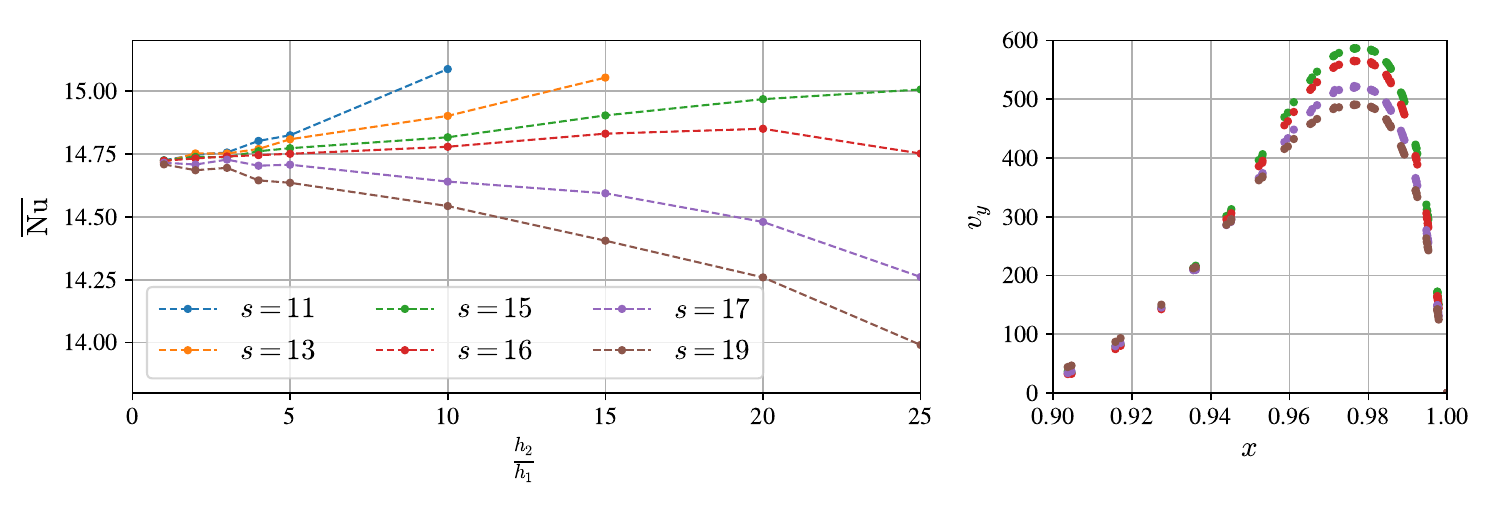}
	\caption{\textit{left:} The average Nusselt number on the cold wall as a function of refinement ratio $\frac{h_2}{h_1}$ for different support sizes $s$. \textit{right:} Vertical velocity cross-section close to the right wall at $y = 0.5$ corresponding to the rightmost cases from the left graph with $\frac{h_2}{h_1} = 25$ and different support sizes $s$. The cross-section shows values in nodes closer than $5h_1$ in $y$ coordinate.
		Support sizes are chosen as $s = k_s N_p + 1$, rounded to the closest integer for $k_s \in \{1.75, 2, 2.25, 2.5, 2.75, 3\}$.}
	\label{fig:refineSupportSize}
\end{figure}

We present an additional observation regarding the numerical method parameterisation that can only be made when dealing with aggressively refined discretisations. The recommended support size $s > 2N_p$ for the approximation, corresponding to $s=12$ in the analysed case, is not sufficient to reach a stable solution for higher refinement ratios $\frac{h_2}{h_1}$ as shown in the left graph of Figure~\ref{fig:refineSupportSize}. The results for the $\mathrm{Ra}=10^5$, $\mathrm{Pr}=100$, $n=0.6$ case calculated with $m=2$ and $h_1 = 0.0025$ show that the range of stable refinement ratios raises with increasing support size until reaching $\frac{h_2}{h_1} = 25$ with $s=15$. Based on this we conservatively choose $s = 2.5 N_p + 1$ rounded to the closest integer, corresponding to $s=16$ in Figure~\ref{fig:refineSupportSize}, as a rule for support size selection in the subsequent refined cases. The Nusselt number for the selected support size remained practically the same for the entire range of refinement ratios resulting in discretisations ranging from 140134 nodes for $\frac{h_2}{h_1} = 1$ to 12438 nodes for $\frac{h_2}{h_1} = 25$.  Note that the $h_\mathrm{max}$ was not enforced for the discretisation and that the dense band width was set to $w=0$, further amplifying the impact of the changing node density on the boundary flow layer. The discrepancy between the different considered support sizes seen in Nusselt numbers on the left is also apparent from the vertical velocity cross-sections shown on the right graph of Figure~\ref{fig:refineSupportSize} with the larger support sizes appearing to smoothen the velocity peak.

\begin{figure}
	\includegraphics[width=\linewidth]{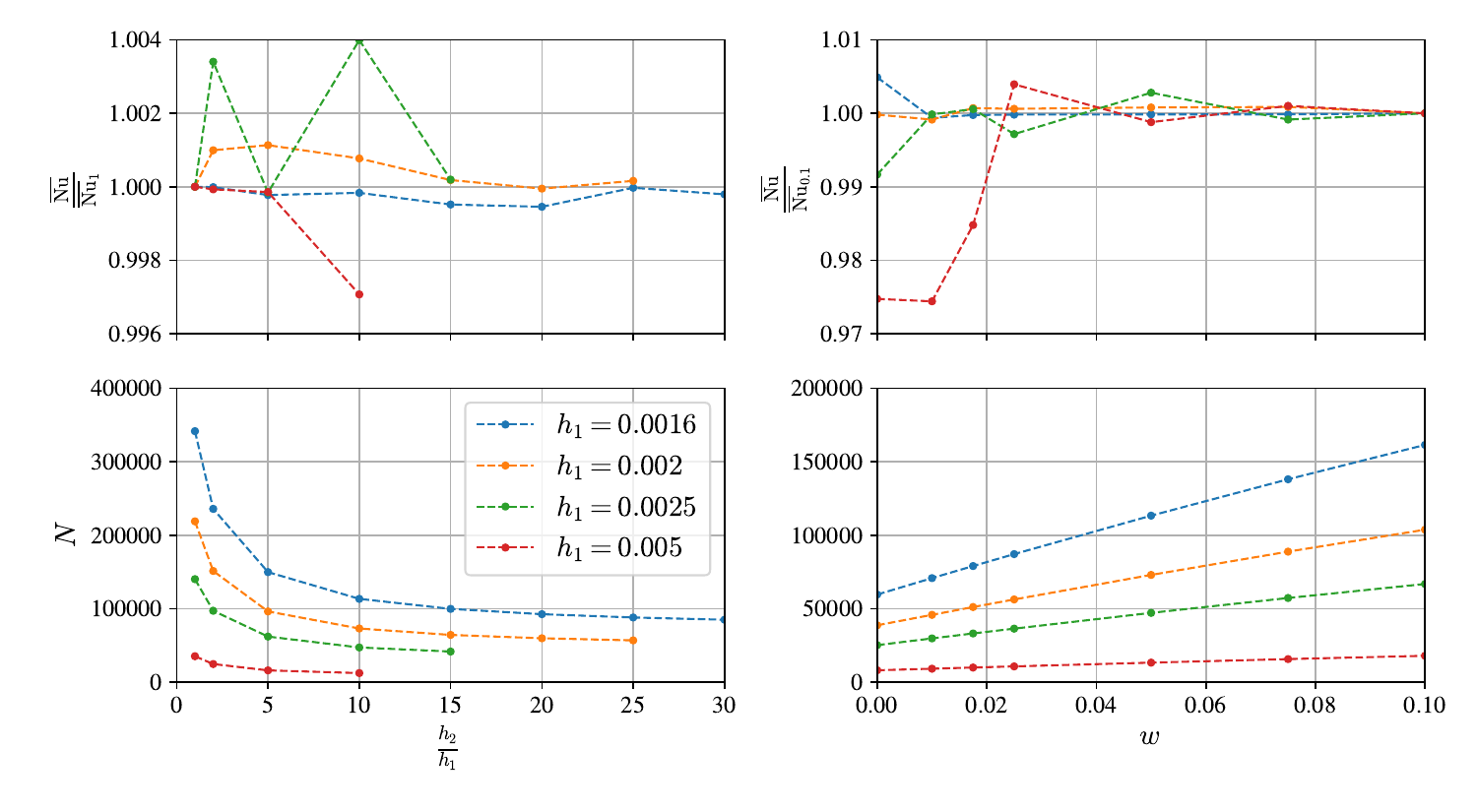}
	\caption{Impact of refinement scheme parameter variation on the Nusselt number show in the top and node count in the bottom graphs. Nusselt values are normed to the value attained with the best considered parameter value. \textit{left:} Variation of the refinement ratio $\frac{h_2}{h_1}$ with the 
		dense boundary band held constant at $w=0.05$. 
		\textit{right:} Variation of 
		the dense boundary band $w$ with the refinement ratio held constant at 
		$\frac{h_2}{h_1} = 10$.}
	\label{fig:refineParameterConvergence}
\end{figure}

The refinement parameters, used in Figure~\ref{fig:refinementScheme} 
and in 
other refined results with unspecified values, were chosen based on a
convergence analysis. The analysis performed at $\mathrm{Ra}=10^6$, $\mathrm{Pr}=100$ 
and $n=0.6$ 
would need to be repeated for other cases with a significant 
differences in 
boundary layer thickness and other flow characteristics. Refinement 
parameters 
have been individually varied for different border densities with 
the results 
shown in Figure~\ref{fig:refineParameterConvergence}. 

We assume that a further refinement in parameters would not lead to significant changes and therefore we norm the Nusselt number, shown in the upper graphs, to the best present value. 
Normalized values allow for comparison between cases with a different maximum density 
and conveniently show the relative offset in value.

The refinement ratio $\frac{h_2}{h_1}$ sweep is performed with a band width $w=0.05$ that is 
wide enough 
to not have a meaningful impact on the results as established in the 
following 
paragraph. The resulting Nusselt numbers, shown in the upper left 
graph of 
Figure~\ref{fig:refineParameterConvergence}, are relatively 
unaffected by the 
refinement ratio, as long as the resulting $h_2$ is
small enough to ensure 
numerical stability. Results are stable, with variations within 1\%. We 
chose $\frac{h_2}{h_1} = 20$ as a conservative choice for refinement aggressiveness, corroborated with Figure~\ref{fig:refineSupportSize} for the selected support size.
The conservative choice of refinement ratio is justified by the graph of computational node count dependence, shown in the lower left graph of Figure~\ref{fig:refineParameterConvergence}, where we can see that further increasing refinement ratio leads to diminishing reduction in the number of computational nodes. This analysis as also used to set the $h_{\mathrm{max}} = 0.05$ based on refinement ratios where stable solution was no longer achieved.

The dense boundary width $w$ sweep is calculated with a refinement ratio of
$\frac{h_2}{h_1} = 10$, that has been determined to be adequate in the previous paragraph, and the results 
are shown in the upper right graph of 
Figure~\ref{fig:refineParameterConvergence}. The main goal of the 
dense 
boundary band is providing a sufficient discretisation for the high 
velocity 
flow and the corresponding gradients in the boundary layer. We can use 
Figure~\ref{fig:logConvergenceProfiles} to estimate the required 
width to 
0.02-0.03, based on the distance from the boundary where flow 
velocity reduces 
to half its maximum value. The estimate is confirmed with numerical 
results 
that show no improvement when increasing the boundary band width $w$ 
beyond a 
point where it covers the boundary layer flow. The effect of $w$ is 
less 
noticeable when using smaller $h_1$ as the linearly decreasing 
density beyond 
the edge of the dense band still provides a sufficient density as 
long as the 
peak of the flow is covered. We chose $w=0.025$ as the refinement 
parameter for 
further use as it provides practically the same results as $w=0.1$, especially 
at smaller 
$h_1$ that we are mainly interested in, while providing a 
significant reduction 
in node count as seen from the lower right graph of 
Figure~\ref{fig:refineParameterConvergence}. 

The symmetric refinement approach is not ideal, but is suitable as a 
proof-of-concept due to the small number of density function 
parameters, which 
simplifies their analysis and selection. A simple improvement would 
be to treat 
the non-insulated and insulated boundaries differently, since the 
latter do not 
exhibit the sharp convective flow layer and can be discretised with 
a lower 
density. The real improvements to refinement can become arbitrarily 
complex 
and focus on deriving an error indicator based either on a previous, 
less 
refined solution or on the properties of the numerical method itself.

\begin{figure}
	\includegraphics[width=\linewidth]{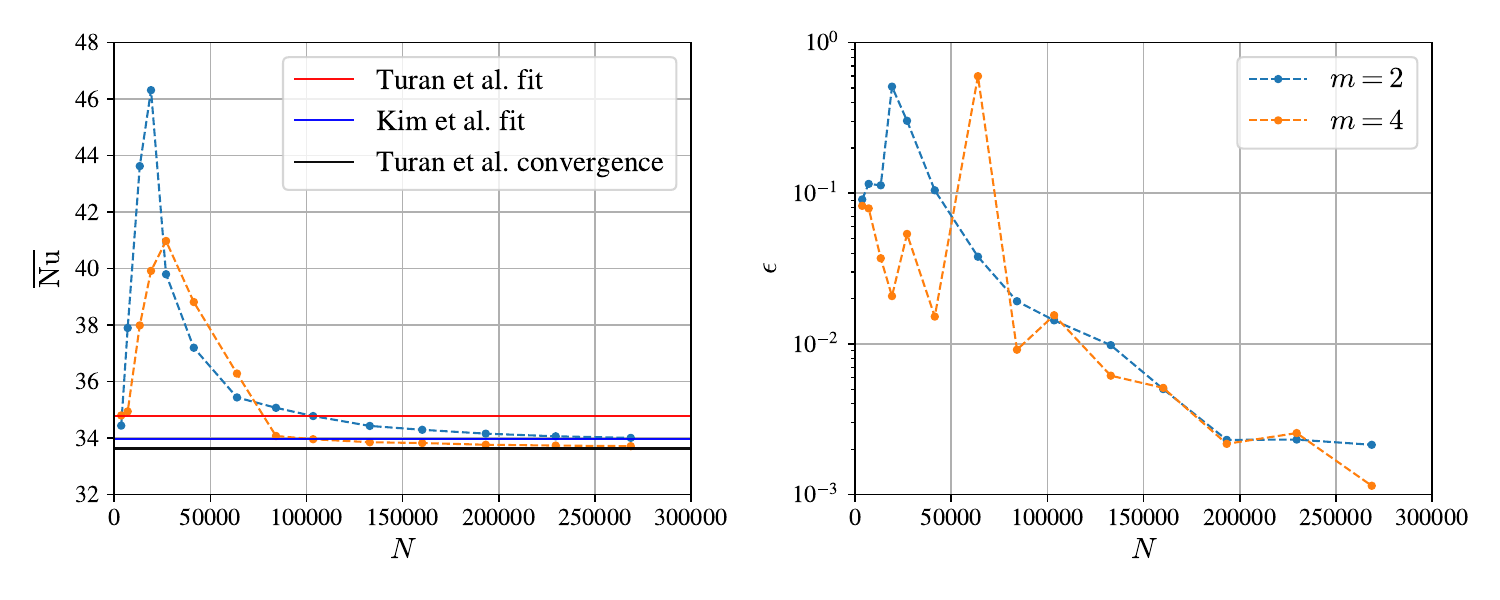}
	\caption{\textit{left:} Convergence study for the problematic 
	$\mathrm{Ra}=10^6$, $\mathrm{Pr}=100$, $n = 0.6$ case repeated with a refined density. Red and blue 
		horizontal lines 
		show Nusset values derived from correlations defined in Eq. 
		\eqref{eq:turanFit} and \eqref{eq:kimFit} while the black 
		shows the best value 
		from the convergence study presented in Turan et al.	
		\protect\cite{turan2011laminar}. \textit{right:} The 
		relative symmetry 
		error for the repeated convergence study.}
	\label{fig:refinedConvergence06}
\end{figure}

We use the refined density to recalculate the convergence for 
$\mathrm{Ra}=10^6$, 
$n=0.6$ case that did not converge with the constant density 
discretisation.  The new convergence results 
are shown in 
Figure~\ref{fig:refinedConvergence06} for the two different considered augmentation orders $m$. Not 
only is convergence achieved, but it is also achieved with a 
significantly 
smaller number of nodes and with a drastically reduced symmetry 
error. The computational times are shown in Table~\ref{tab:timing} to further emphasise that the reduction in the number of computational nodes translates into savings in computational time required. Even with the increased support size the refined discretisation still provides a solution in more than 10 times shorter time that the unrefined. 
The timing was performed on 4 cores of \texttt{Intel(R) Xeon(R) E5520} CPU with the frequency fixed to 2.27GHz.

\begin{table}[h]
	\begin{center}
		\begin{tabular}{|c|c|c|c|c|c|c|c|}
			\hline
			                           & $h_1$                   & $h_2$                   & $N$                     & $m$ & $s$ & $\overline{\mathrm{Nu}}$ & $t$[h]  \\
			\hline
			\multirow{2}{*}{refined}   & \multirow{4}{*}{0.0025} & \multirow{2}{*}{0.05}   & \multirow{2}{*}{26961}  & 2 & 16 & 39.8 & \bf{5} \\
			\cline{5-8}
			                           &                         &                         &                         & 4 & 39 & 41.0 & \bf{9} \\
			\cline{1-1} \cline{3-8}
			\multirow{2}{*}{unrefined} &                         & \multirow{2}{*}{0.0025} & \multirow{2}{*}{140134} & 2 & 13 & 38.8 & \bf{63} \\
			\cline{5-8}
			                           &                         &                         &                         & 4 & 31 & 41.0 & \bf{182} \\
			\hline
		\end{tabular}
	\end{center}
	\caption{A table showing differences in computational time between the refined and the unrefined discretisation and between the different approximation orders for the $\mathrm{Ra}=10^6$, $\mathrm{Pr}=100$, $n = 0.6$ case.}
	\label{tab:timing}
\end{table}

\subsection{Comparison with reference data}
Non-Newtonian fluid dynamics in a differentially heated cavity have 
already 
been tackled previously, and we have performed our study on a 
matching case 
for verification. We use the results provided by Turan et 
al. \cite{turan2011laminar} and Kim et al. \cite{kim2003transient}, 
who used 
FVM with SIMPLE coupling and Upwind or QUICK stabilisation. We 
compare those 
results with a refined RBF-FD and ACM, without 
stabilisation of 
convective terms to minimise the effect of numerical diffusion.

Both publications provide an empirical fit for the average Nusselt 
number as a 
function of Rayleigh number Ra, Prandtl number Pr and non-Newtonian 
exponent 
$n$
\begin{align}
	\overline{\mathrm{Nu}}_\mathrm{Turan} &= 0.162 
	\mathrm{Ra}^{0.043} 
	\frac{\mathrm{Pr}^{0.341}}{(1 + \mathrm{Pr})^{0.091}} \left( 
	\frac{\mathrm{Ra}^{(2-n)}}{\mathrm{Pr}^n} 
	\right)^\frac{1}{2(n+1)} 
	\exp(\mathrm{C}(n-1)), \\ \label{eq:turanFit}
	\mathrm{C} &= \left\{\begin{array}{ll}
		1.343 \mathrm{Ra}^{0.065} \mathrm{Pr}^{0.036} & \text{where 
		} n \le 1 \\
		0.858 \mathrm{Ra}^{0.071} \mathrm{Pr}^{0.034} & 
		\mathrm{otherwise}
	\end{array}\right. ,\\
	\overline{\mathrm{Nu}}_\mathrm{Kim} &= 0.3 n^{0.4} 
	\mathrm{Ra}^\frac{1}{3n 
		+ 1}, \label{eq:kimFit}
\end{align}
that we can use to compare against our results.

The Nusselt values calculated by the provided correlations are added 
to 
Figure~\ref{fig:refinedConvergence06}, supplemented by the exact 
numerical 
value from the convergence study performed by Turan et al. on this 
case. All 
results agree well, with small deviations that are normal due to 
discretisation 
errors and different numerical methods.
We must also keep in mind 
that the  
provided correlations, in all their complexity, are still only 
empirical fits 
over a wide range of flow regimes and are a rather crude approximation, as can be 
seen from 
the relatively large discrepancy between the Turan et al. fit and the convergence 
value from the 
corresponding study.

\begin{figure}
	\includegraphics[width=\linewidth]{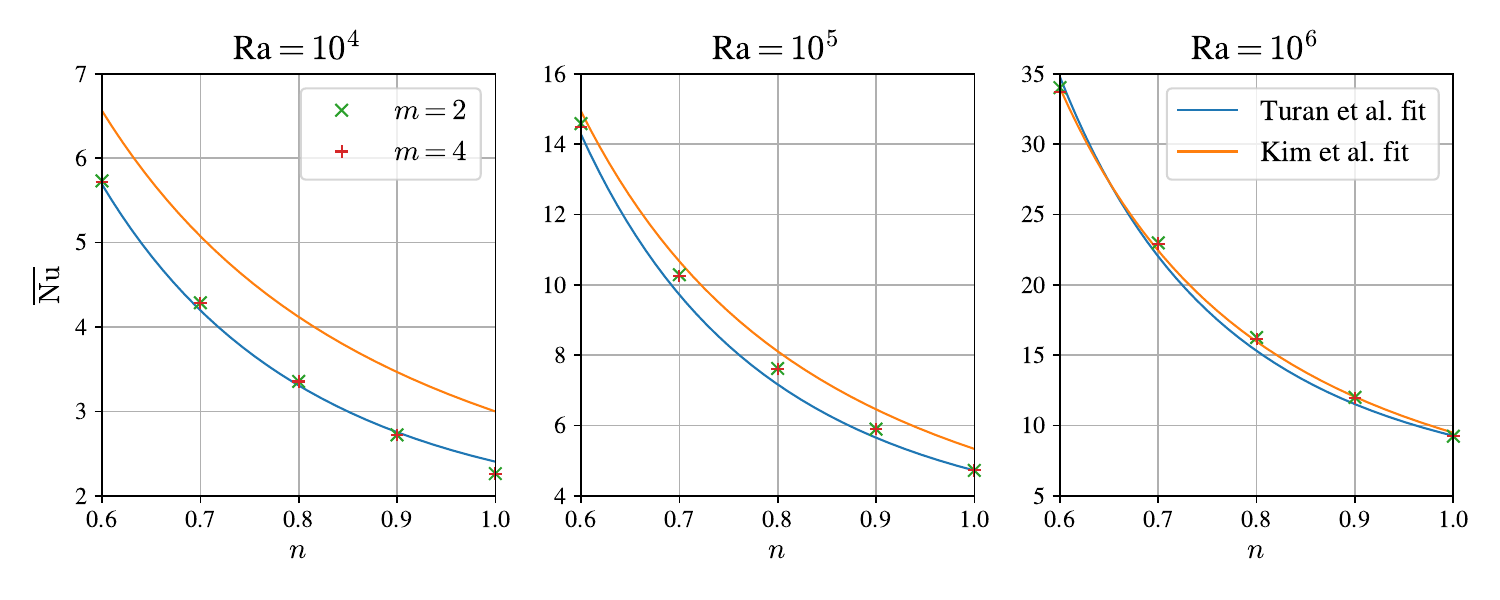}
	\caption{Comparison between the calculated Nusselt values and 
	correlations 
		provided in existing literature. Solid lines represent 
		correlations defined 
		in equations \eqref{eq:turanFit} and \ref{eq:kimFit}.}
	\label{fig:fitComparison}
\end{figure}

We extend the fit comparison to all cases considered with the results shown in Figure~\ref{fig:fitComparison}. All computed values fall within the range of the provided correlation functions, again confirming that the meshless RBF-FD method yields valid results. The relative difference between the fits increases significantly for calmer cases with a smaller Nusselt number, casting further 
doubt on the accuracy of the empirical fits for a wide range of parameters. The case-dependent discrepancy between our numerical results and the provided fits is similar to that in the source publication.

We also tested different monomial orders for the entire range of considered cases shown in Figure~\ref{fig:fitComparison}. The difference between the results computed with augmentation orders is in all cases smaller than the difference between Turan et al. and Kim et al., which means that the lower order $m=2$ is adequate for a stable and accurate computation when the discretisation is sufficiently dense. Whether using higher order method is beneficial is case dependent as can be seen from Figure~\ref{fig:nusseltConvergence}. Cases with lower Ra, where resolution of the boundary layer is not problematic, appear to benefit from the higher order by reaching node density independent result at a drastically lower node count. The opposite is true for the high Ra and low $n$ cases where increasing the order is at times counterproductive. Using a higher order method appears beneficial for the refined discretisation cases shown in Figure~\ref{fig:refinedConvergence06} but the final benefits of increasing the augmentation order are again case dependent and require additional analysis to determine optimal compromise between the lower total number of computational nodes and the increased computational time required per computational node.

We report the average Nusselt value at the cold boundary calculated with the finest discretisation for the full range of considered parameters at $\mathrm{Pr}=100$ in Table~\ref{tab:results}.

\begin{table}[h]
	\begin{center}
		\begin{tabular}{|c|c|c|c|c|c|}
			\hline
			Ra & $n = 0.6$ & $n = 0.7$ & $n = 0.8$ & $n = 0.9$ & $n = 1$ \\
			\hline
			$10^{4}$ & 5.72 & 4.28 & 3.35 & 2.72 & 2.26 \\
			\hline
			$10^{5}$ & 14.49 & 10.26 & 7.61 & 5.90 & 4.72 \\
			\hline
			$10^{6}$ & 33.71 & 22.92 & 16.16 & 11.97 & 9.23 \\
			\hline
		\end{tabular}
	\end{center}
	\caption{A table of the best obtained average Nusselt values for 
		all of the 
		considered parameters. Presented values are calculated with 
		$m=2$. Refined 
		density discretisation is used for the problematic $\mathrm{Ra}=10^6$, $n \in \{0.6, 
		0.7\}$ cases.}
	\label{tab:results}
\end{table}

\section{Geometrical flexibility}
\label{ch:flexibility}
\begin{figure}
	\includegraphics[width=\linewidth]{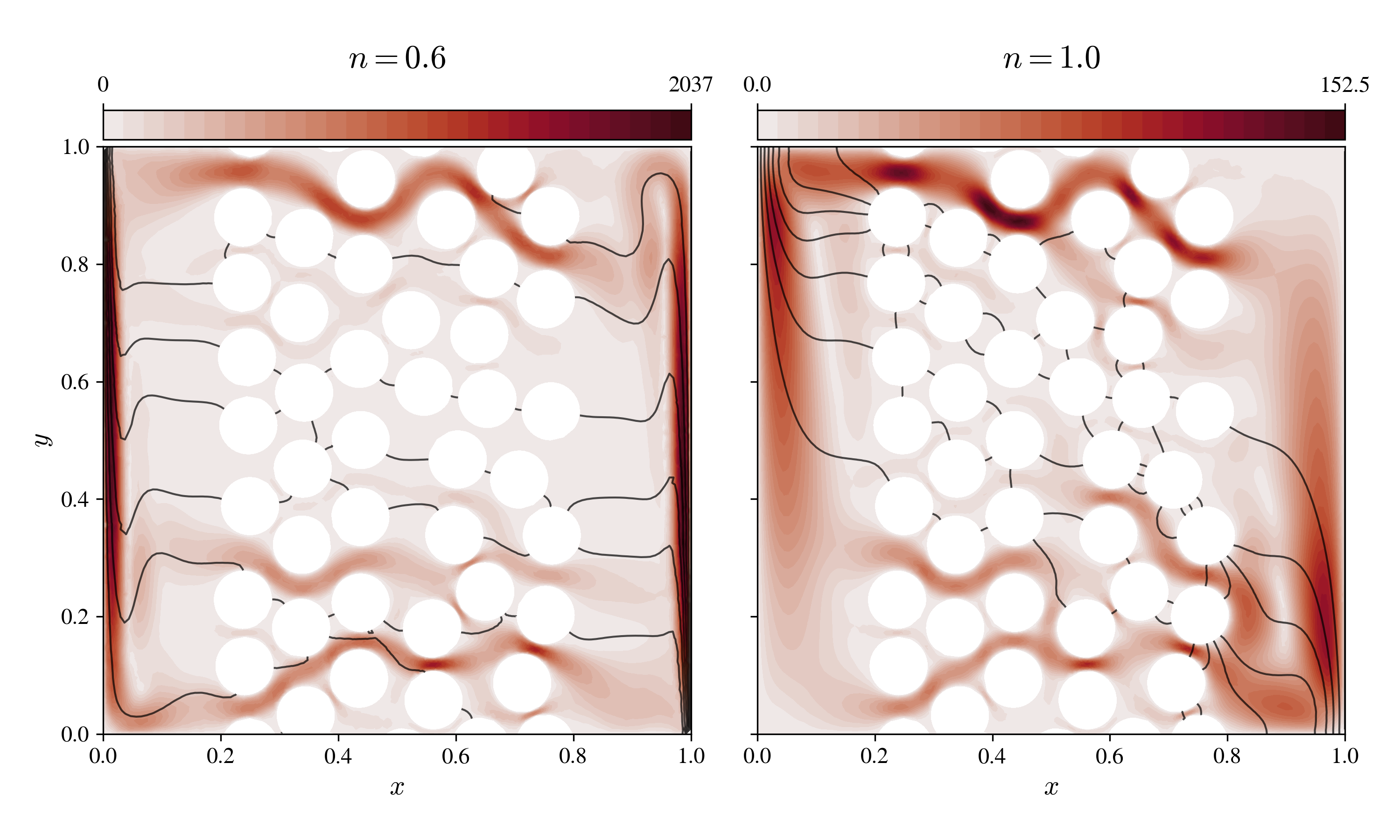}
	\caption{Flow profiles for $\mathrm{Ra}=10^6$, $\mathrm{Pr}=100$ case with a 
	porous central 
		section. Velocity magnitude is visualised with a heat-map 
		while the 
		overlaid contours display the changes in temperature. The 
		strongest 
		shear-thinning $n=0.6$ case is shown on the left sub-figure 
		while the 
		right 
		displays the Newtonian $n=1$ case.} 
	\label{fig:porous2D}
\end{figure}

Finally, we apply the described solution procedure to  more complex 
domains in order to demonstrate the geometrical flexibility that we 
touted as 
one of the main benefits. In the first example, shown in 
Figure \ref{fig:porous2D} we add circular obstructions to the 
central part of 
the De Vahl Davis case to simulate how the non-Newtonian behaviour 
would impact 
the convective flow passing through a porous filter. The 
computational node 
distribution used in this case, shown in the right panel of 
Figure~\ref{fig:nodePlacement}, utilises a modified refinement 
strategy that 
increases node density in narrow channels and on nearby boundaries 
to ensure 
that at least two computational nodes discretise the channel's 
width. Even 
though the dimensionless numbers defined in section 
\ref{ch:problemFormulation} 
are no longer suitable for this case we stick with them to enable 
comparison 
with unobstructed flow profiles shown in Figure~\ref{fig:flowFields}. We use the 
same random 
filter configuration for the shear-thinning non-Newtonian fluid on 
the left and 
the Newtonian fluid on the right graph of Figure~\ref{fig:porous2D}. The results are unsurprising with both the change in temperature profile and the reduction of maximum velocity confirming that the filter's narrow channels present a far greater hindrance for the Newtonian fluid. This can further be confirmed by comparing the average Nusselt number on the cold wall against the unobstructed case with the same $h=0.008$ discretisation. The $\overline{\mathrm{Nu}}$ value decreases in both cases, indicating weaker convective heat transfer, but the change in the non-Newtonian $n=0.6$ case from 35.8 to 35 is almost negligible while the $\overline{\mathrm{Nu}}$ value practically halves from 9.5 to 4.8 in the Newtonian $n=1$ case.

\begin{figure}
	\centering
	\includegraphics[width=0.8\linewidth]{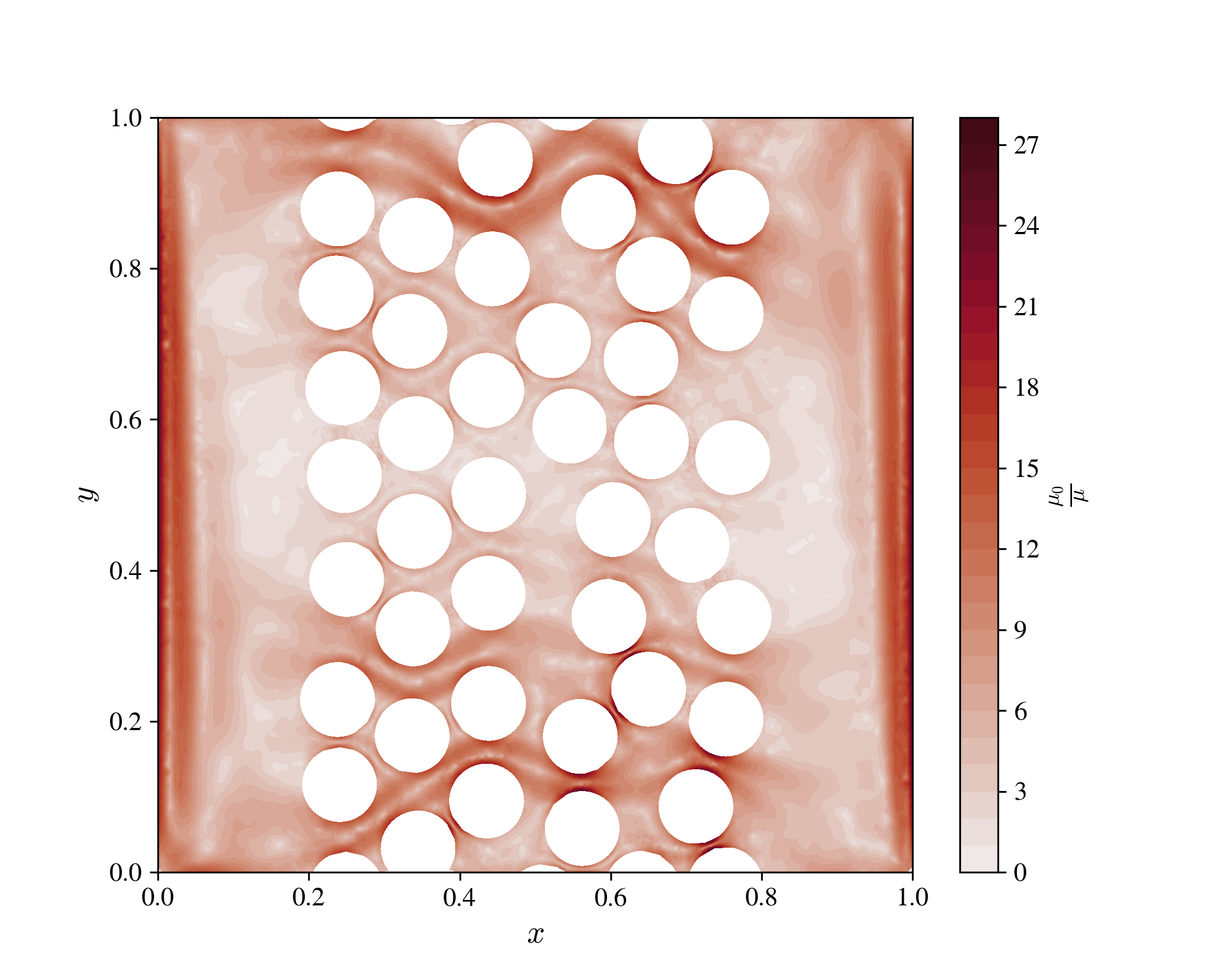}
	\caption{Variation in viscosity due to shear-thinning effects. 
	Higher 
		values signify lower viscosity due to the normed inverse 
		visualisation 
		that is chosen to highlight the areas where shear rate 
		caused thinning 
		is the strongest.}
	\label{fig:porous2Dviscosity}
\end{figure}

The shear-thinning filter case presents a suitable opportunity to 
examine the 
spatial variation in viscosity shown in 
Figure~\ref{fig:porous2Dviscosity}. We 
show the inverse of viscosity to better highlight the shear-thinning 
aspect 
with the most affected areas corresponding to edges of the high 
velocity layer 
next to the vertical walls and in high velocity channels through the 
filter as 
seen in Figure~\ref{fig:porous2D}. Viscosity behaviour next to the 
vertical 
walls is identical to the non filter cases.

\begin{figure}
	\includegraphics[width=\linewidth]{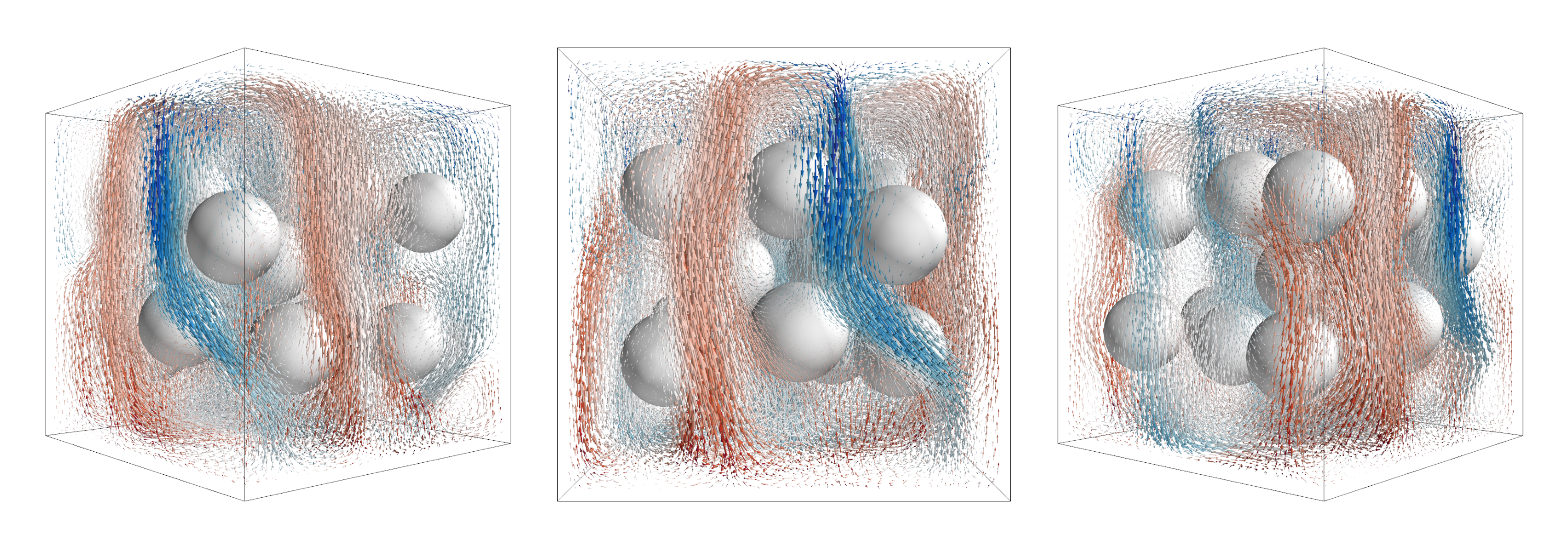}
	\caption{Natural convection in a 3D differentially heated cavity with obstructions. This case differs from the previous as the 
		temperature differential and resulting flow is now vertical. The glyphs are sized according to the velocity magnitude in the computational node and coloured according to the temperature complying with the colourbar in Figure~\ref{fig:porous3Dintersections}.}
	\label{fig:porous3D}
\end{figure}

\begin{figure}
	\includegraphics[width=\linewidth]{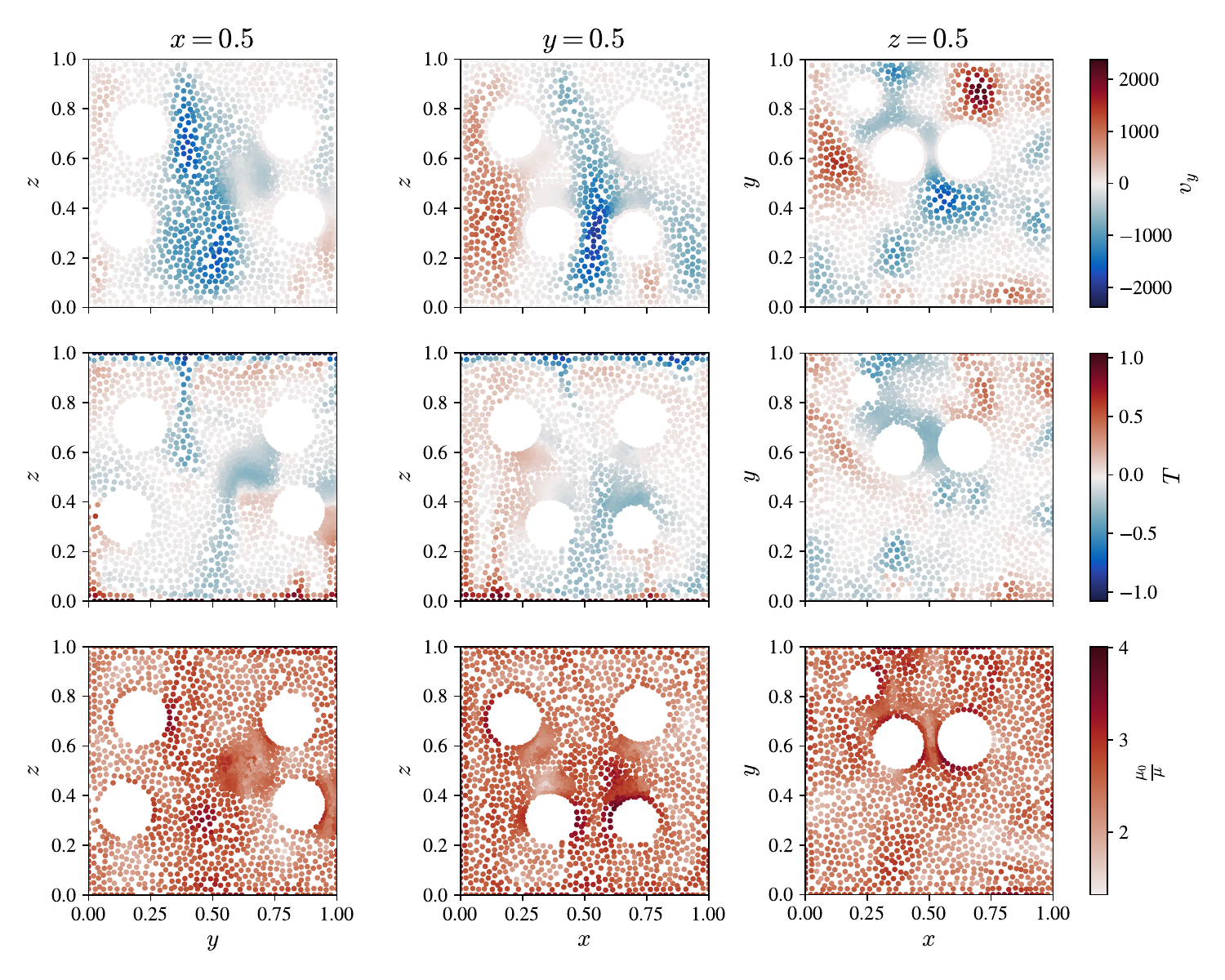}
	\caption{Plane intersections of the 3D case shown in 
		Figure~\ref{fig:porous3D}. The first row displays vertical 
		velocity, the 
		second row the temperature and the third row the inverse of 
		the 
		shear-thinning viscosity. Values and positions shown on the 
		scatter plot 
		correspond to computational nodes within $\frac{h}{2}$ of 
		the intersecting 
		plane.}
	\label{fig:porous3Dintersections}
\end{figure}

We formulated the method in a dimensionally agnostic manner, which 
can be with 
some effort directly transformed into dimensionally agnostic 
object-oriented 
code design~\cite{medusa}, using powerful C++ template systems. The 
filtered 2D 
case can thus be quickly transformed into its 3D variant shown in 
Figure~\ref{fig:porous3D} with minimal interference  in the core 
code that is 
limited only to the parameters (refer to the supplied repository for 
the actual 
code). For better visual representation plane intersections of the 
vertical 
velocity, temperature and viscosity are presented in 
Figure~\ref{fig:porous3Dintersections}. The purpose of presenting 
this 3D 
results is twofold. First, we want to demonstrate the ability of the 
presented 
method to address a complex physical problem on an irregular 3D 
domain. Second, 
we want to demonstrate the powerful concept of merging the generally 
formulated 
RBF-FD method with a generic programming to design an dimension 
independent 
solution procedure. More details on the C++ implementation aspect 
can be found 
in~\cite{medusa}.  

\section{Conclusions}
In this paper, we have proposed a dimension independent refined meshless solution procedure for NNC. The main advantage of the proposed approach is that it can operate on scattered nodes, which greatly facilitates the consideration of complex 3D domains and the implementation of refined discretisation, both of which are demonstrated in this paper. Moreover, the approximation weights are computed individually, allowing for a variation in the stencil size $s$, monomial augmentation order $m$, the type and order of the RBFs used, etc. The ease of adaptation on all mentioned levels allows the use of a sophisticated but slow approximation only where it is truly needed with faster alternatives elsewhere, which is often the majority of the domain when dealing with realistic cases.

From the meshless method point of view, we have shown that the method is appropriate for the considered problem given that the discretisation is sufficiently dense to describe the present flow structures. We have shown how the problems with inadequate discretisation manifest themselves and how they can be resolved by refinement. We have shown that the proposed solution procedure allows for an aggressively refined discretisation that can be pushed even further if the support size is slightly increased from the usual recommendation, leading to an order of magnitude lower computational times. 

From the non-Newtonian fluid dynamics point of view, we presented new results for the benchmark case that is already solved with two similar FVM approaches~\cite{turan2011laminar,kim2003transient} with a fundamentally different numerical approach and, most importantly, without stabilisation, that effectively introduces numerical diffusion in the solution. 

In future work, we will expand the analysis of how the interplay between node density and stencil size influences sharp flow structures, like the boundary layer shown in this paper, and attempt to determine what are the requirements for an accurate reproduction. We will use this knowledge to develop an h-adaptive solution, which requires an adequate error indicator and a refinement logic that constructs the target node density function from the error indicator data and refines the nodes accordingly~\cite{slak2019adaptive}.

In addition, we also plan to analyse how hyperviscosity~\cite{shankar2018hyperviscosity} and adaptive upwind~\cite{kosec2014simulation} stabilisations, and using different approximation approaches~\cite{hativc2021meshless} affect the flow structures.

\begin{acknowledgments}
	The authors would like to acknowledge the financial support of 
	the Slovenian Research and Innovation Agency (ARIS) research core funding No.\ 
	P2-0095, Young Researcher programme PR-10468, and project funding N2-0275. Funded by National Science Centre, Poland under the OPUS call in the Weave programme 2021/43/I/ST3/00228.
	This research was funded in whole or in part by National Science Centre (2021/43/I/ST3/00228). For the purpose of Open Access,
	the author has applied a CC-BY public copyright licence to any Author Accepted Manuscript (AAM) version arising from this submission.
\end{acknowledgments}

\section*{Data Availability Statement}

The code for the presented numerical methodology, plotting scripts, and data that support the findings of this study are openly 
available on Zenodo at http://doi.org/10.5281/zenodo.14901942, reference number~\cite{rot2025zenodo}. For ease of access, the code and data are also available in the GitLab repository, reference number \cite{git}.

\bibliographystyle{aip}
\bibliography{RNNC}

\begin{thebibliography}{10}

\bibitem{rahimi2018comprehensive}
A.~Rahimi, A.~D. Saee, A.~Kasaeipoor, and E.~H. Malekshah,
\newblock International Journal of Numerical Methods for Heat \& Fluid Flow
  (2018).

\bibitem{bejan2013convection}
A.~Bejan,
\newblock {\em Convection heat transfer},
\newblock John Wiley \& Sons, 2013.

\bibitem{Chhabra2010}
R.~P. Chhabra,
\newblock {\em Non-Newtonian Fluids: An Introduction}, pages 3--34,
\newblock Springer New York, New York, NY, 2010.

\bibitem{Wang2010}
X.~Wang and X.~Li,
\newblock Finite Elements in Analysis and Design {\bf 46}, 551  (2010).

\bibitem{Charm1965}
S.~Charm and G.~Kurland,
\newblock Nature {\bf 206}, 617 (1965).

\bibitem{Gratao2007}
A.~Gratão, V.~Silveira, and J.~Telis-Romero,
\newblock Journal of Food Engineering {\bf 78}, 1343  (2007).

\bibitem{WeltiChanes2005}
J.~Welti-Chanes, F.~Vergara-Balderas, and D.~Bermúdez-Aguirre,
\newblock Journal of Food Engineering {\bf 67}, 113  (2005),
\newblock IV Iberoamerican Congress of Food Engineering (CIBIA IV).

\bibitem{bingham1917investigation}
E.~C. Bingham,
\newblock {\em An investigation of the laws of plastic flow},
\newblock Number 278, US Government Printing Office, 1917.

\bibitem{yang2019comprehensive}
L.~Yang and K.~Du,
\newblock Journal of Thermal Analysis and Calorimetry , 1 (2019).

\bibitem{Kwack2014}
J.~Kwack and A.~Masud,
\newblock Computational Mechanics {\bf 53}, 751 (2014).

\bibitem{patankar1980numerical}
S.~Patankar,
\newblock Numerical heat transfer and fluid flow, taylor \& francis london,
  1980.

\bibitem{kim2003transient}
G.~B. Kim, J.~M. Hyun, and H.~S. Kwak,
\newblock International Journal of Heat and Mass Transfer {\bf 46}, 3605
  (2003).

\bibitem{turan2011laminar}
O.~Turan, A.~Sachdeva, N.~Chakraborty, and R.~J. Poole,
\newblock Journal of Non-Newtonian Fluid Mechanics {\bf 166}, 1049 (2011).

\bibitem{bozorg2019two}
M.~V. Bozorg and M.~Siavashi,
\newblock International Journal of Mechanical Sciences {\bf 151}, 842 (2019).

\bibitem{moraga2017geometric}
N.~O. Moraga, M.~A. Marambio, and R.~C. Cabrales,
\newblock International Communications in Heat and Mass Transfer {\bf 88}, 108
  (2017).

\bibitem{vasco2014parallel}
D.~A. Vasco, N.~O. Moraga, and G.~Haase,
\newblock Numerical Heat Transfer, Part A: Applications {\bf 66}, 990 (2014).

\bibitem{loenko2019natural}
D.~S. Loenko, A.~Shenoy, and M.~A. Sheremet,
\newblock Energies {\bf 12}, 2149 (2019).

\bibitem{alsabery2017transient}
A.~Alsabery, A.~Chamkha, H.~Saleh, and I.~Hashim,
\newblock Powder Technology {\bf 308}, 214 (2017).

\bibitem{mishra2018natural}
L.~Mishra and R.~P. Chhabra,
\newblock Heat Transfer Engineering {\bf 39}, 819 (2018).

\bibitem{kefayati2014simulation}
G.~R. Kefayati,
\newblock International Communications in Heat and Mass Transfer {\bf 53}, 139
  (2014).

\bibitem{liu2002mesh}
G.-R. Liu,
\newblock {\em Mesh free methods: moving beyond the finite element method},
\newblock CRC press, 2002.

\bibitem{tolstykh2003using}
A.~I. Tolstykh and D.~A. Shirobokov,
\newblock Computational Mechanics {\bf 33}, 68 (2003).

\bibitem{slak2019generation}
J.~Slak and G.~Kosec,
\newblock SIAM Journal on Scientific Computing {\bf 41}, A3202 (2019).

\bibitem{trojak2022acm}
W.~Trojak, N.~Vadlamani, J.~Tyacke, F.~Witherden, and A.~Jameson,
\newblock Computers \& Fluids {\bf 247}, 105634 (2022).

\bibitem{kosec2018localIrregularFlow}
G.~Kosec,
\newblock Advances in Engineering Software {\bf 120}, 36 (2018),
\newblock Civil-Comp - Part 2.

\bibitem{yasuda2023bulkViscosityACM}
T.~Yasuda, I.~Tanno, T.~Hashimoto, K.~Morinishi, and N.~Satofuka,
\newblock Computers \& Fluids {\bf 258}, 105885 (2023).

\bibitem{slak2019adaptive}
J.~Slak and G.~Kosec,
\newblock International Journal for Numerical Methods in Engineering {\bf 119},
  661 (2019).

\bibitem{jancic2021p}
M.~Jan{\v{c}}i{\v{c}}, J.~Slak, and G.~Kosec,
\newblock p-refined rbf-fd solution of a poisson problem,
\newblock in {\em 2021 6th International Conference on Smart and Sustainable
  Technologies (SpliTech)}, pages 01--06, IEEE, 2021.

\bibitem{anderson2010fundamentals}
J.~D. Anderson~Jr,
\newblock {\em Fundamentals of aerodynamics},
\newblock Tata McGraw-Hill Education, 2010.

\bibitem{Tritton1988}
D.~J. Tritton,
\newblock {\em Physical Fluid Dynamics},
\newblock Oxford Science Publ, Clarendon Press, 1988.

\bibitem{Note1}
The shear rate tensor for incompressible fluids is the same as the strain rate
  tensor, which can be expressed as the symmetric part of the velocity
  gradient.

\bibitem{Note2}
The limit is arbitrarily chosen and has a negligible effect on the result \cite
  {kim2003transient}.

\bibitem{DeVahlDavis1983}
G.~De~Vahl~Davis,
\newblock International Journal for Numerical Methods in Fluids {\bf 3}, 249
  (1983).

\bibitem{hardin2004discretizing}
D.~P. Hardin and E.~B. Saff,
\newblock Notices of the AMS {\bf 51}, 1186 (2004).

\bibitem{shankar2018robust}
V.~Shankar, R.~M. Kirby, and A.~L. Fogelson,
\newblock SIAM J. Sci. Comput. {\bf 40}, 2584 (2018).

\bibitem{depolli2022parallel}
M.~Depolli, J.~Slak, and G.~Kosec,
\newblock Parallel domain discretization algorithm for rbf-fd and other
  meshless numerical methods for solving pdes, 2022.

\bibitem{duh2020fast}
U.~Duh, G.~Kosec, and J.~Slak,
\newblock {SIAM} Journal on Scientific Computing {\bf 43}, A980 (2021).

\bibitem{duh_discretization_2024}
U.~Duh, V.~Shankar, and G.~Kosec,
\newblock Journal of Scientific Computing {\bf 100}, 51 (2024).

\bibitem{davydov2023stencilSelection}
O.~Davydov, D.~T. Oanh, and N.~M. Tuong,
\newblock Journal of Computational and Applied Mathematics {\bf 425}, 115031
  (2023).

\bibitem{flyer2016polynomials1}
N.~Flyer, B.~Fornberg, V.~Bayona, and G.~A. Barnett,
\newblock Journal of Computational Physics {\bf 321}, 21 (2016).

\bibitem{Note3}
The 6 monomials in 2-D case with $m=2$ would be $q = \{1, x, y, x^2, xy,
  y^2\}$.

\bibitem{jancic2021monomial}
M.~Jan{\v{c}}i{\v{c}}, J.~Slak, and G.~Kosec,
\newblock Journal of Scientific Computing {\bf 87}, 9 (2021).

\bibitem{leborne2023guidelinesRBFFD}
S.~Le~Borne and W.~Leinen,
\newblock Journal of Scientific Computing {\bf 95}, 8 (2023).

\bibitem{bayona2017}
V.~Bayona, N.~Flyer, B.~Fornberg, and G.~A. Barnett,
\newblock Journal of Computational Physics {\bf 332}, 257 (2017).

\bibitem{jancic_strong_2023}
M.~Jančič and G.~Kosec,
\newblock Engineering with Computers  (2023).

\bibitem{chorin1967acm}
A.~J. Chorin,
\newblock Journal of Computational Physics {\bf 2}, 12 (1967).

\bibitem{kosec_super_2014}
G.~Kosec, M.~Depolli, A.~Rashkovska, and R.~Trobec,
\newblock Computers \& Structures {\bf 133}, 30 (2014).

\bibitem{kajzer2018GPU}
A.~Kajzer and J.~Pozorski,
\newblock Computers \& Mathematics with Applications {\bf 76}, 997 (2018).

\bibitem{rahman2008ACM}
M.~M. Rahman and T.~Siikonen,
\newblock International Journal for Numerical Methods in Engineering {\bf 75},
  1320 (2008).

\bibitem{medusa}
J.~Slak and G.~Kosec,
\newblock ACM Trans. Math. Softw. {\bf 47} (2021).

\bibitem{shankar2018hyperviscosity}
V.~Shankar and A.~L. Fogelson,
\newblock Journal of Computational Physics {\bf 372}, 616  (2018).

\bibitem{kosec2014simulation}
G.~Kosec and B.~{\v{S}}arler,
\newblock Engineering analysis with boundary elements {\bf 45}, 36 (2014).

\bibitem{hativc2021meshless}
V.~Hati{\v{c}}, B.~Mavri{\v{c}}, and B.~{\v{S}}arler,
\newblock Engineering Analysis with Boundary Elements {\bf 131}, 86 (2021).

\bibitem{rot2025zenodo}
M.~Rot and G.~Kosec,
\newblock Refined radial basis function-generated finite difference analysis of
  non-newtonian natural convection, 2025.

\bibitem{git}
Gitlab repository,
\newblock
  \url{https://gitlab.com/e62Lab/public/2022_p_refinednonnewtonianconvection}.

\end{thebibliography}

\end{document}